\documentclass[reqno]{amsproc}
\usepackage{bbm}
\usepackage{amssymb,hyperref,color}

\newtheorem{theorem}{Theorem}[section]
\newtheorem{lemma}[theorem]{Lemma}
\newtheorem{corollary}[theorem]{Corollary}

\theoremstyle{remark}

\numberwithin{equation}{section}

\allowdisplaybreaks
\title[On multiplier analogues of $C+H^\infty$]
{On multiplier analogues of the algebra $C+H^\infty$
on\\ weighted rearrangement-invariant sequence spaces}

\dedicatory{To Professor Yuri Brudnyi on the occasion of his 90th birthday}

\author[O. Karlovych]{Oleksiy Karlovych}
\address{
Centro de Matem\'atica e Aplica\c{c}\~oes\\
Departamento de Matem\'atica\\
Faculdade de Ci\^encias e Tecnologia\\
Universidade Nova de Lisboa\\
Quinta da Torre\\
2829--516 Caparica\\
Portugal} \email{oyk@fct.unl.pt}

\author[S. M. Thampi]{Sandra Mary Thampi}
\address{Centro de Matem\'atica e Aplica\c{c}\~oes,\\
Departamento de Matem\'a\-tica,\\
Faculdade de Ci\^encias e Tecnologia,\\
Universidade Nova de Lisboa,\\
Quinta da Torre,\\
2829--516 Caparica, Portugal}
\email{s.thampi@campus.fct.unl.pt}
\begin{document}
\begin{abstract}
Let $X(\mathbb{Z})$ be a reflexive rearrangement-invariant Banach sequence
space with nontrivial Boyd indices $\alpha_X,\beta_X$ and let $w$ be a 
symmetric weight in the intersection of the Muckenhoupt classes 
$A_{1/\alpha_X}(\mathbb{Z})$ and $A_{1/\beta_X}(\mathbb{Z})$. Let
$M_{X(\mathbb{Z},w)}$ denote the collection of all periodic distributions
$a$ generating bounded Laurent operators $L(a)$ on the space
$X(\mathbb{Z},w)=\{\varphi:\mathbb{Z}\to\mathbb{C}:\varphi w\in X(\mathbb{Z})\}$.
We show that $M_{X(\mathbb{Z},w)}$ is a Banach algebra. Further, we
consider the closure of trigonometric polynomials in $M_{X(\mathbb{Z},w)}$ 
denoted by  $C_{X(\mathbb{Z},w)}$ and $H_{X(\mathbb{Z},w)}^{\infty,\pm}=
\{a\in M_{X(\mathbb{Z},w)}:\widehat{a}(\pm n)=0 \mbox{ for }n<0\}$.
We prove that
$C_{X(\mathbb{Z},w)}+H_{X(\mathbb{Z},w)}^{\infty,\pm}$ are closed
subalgebras of $M_{X(\mathbb{Z},w)}$.
\end{abstract}
\maketitle
\section{Introduction}
Let $C$ and $L^\infty$ be the Banach algebras of {continuous} 
and essentially bounded functions on the unit circle, respectively, and
\[
H^\infty :=\left\{f\in L^\infty\ :\ \widehat{f}(n)=0
\ \mbox{for}\ n<0\right\}.
\]
In 1967, Sarason observed that $C+H^\infty$ is a closed subalgebra of
$L^\infty$ (see \cite[p.~191]{S67} and also \cite[p.~290]{S73}).
An elegant alternative proof of this fact was given by Zalcman 
\cite[Theorem~6.1]{Z69}. A few years later, Rudin \cite[Theorem~1.2]{R75}
presented the following abstract version of Zalcman's argument
(see also \cite[Ch.~VII, Section~3, p.~149]{K98} and 
\cite[Lemma~2.52]{BS06}).

In the theorem below and in what follows, for a Banach space $\mathcal{X}$,
let $\mathcal{B}(\mathcal{X})$ denote the Banach algebra of all bounded
linear operators on $\mathcal{X}$ equipped with the operator norm.
\begin{theorem}\label{th:Zalcman-Rudin}
Suppose $\mathcal{Y}$ and $\mathcal{Z}$ are closed subspaces of a Banach space
$\mathcal{X}$, and suppose there is a collection 
$\Phi\subset\mathcal{B}(\mathcal{X})$ with the following properties:
\begin{enumerate}
\item[(i)]
Every $\Lambda\in\Phi$ maps $\mathcal{X}$ into $\mathcal{Y}$.

\item[(ii)]
Every $\Lambda\in\Phi$ maps $\mathcal{Z}$ into $\mathcal{Z}$.

\item[(iii)]
$\displaystyle \sup_{\Lambda\in\Phi}
\|\Lambda\|_{\mathcal{B}(\mathcal{X})}<\infty$.

\item[(iv)]
To every $y\in\mathcal{Y}$ and to every $\varepsilon>0$ corresponds a
$\Lambda\in\Phi$ such that $\|y-\Lambda y\|_{\mathcal{X}}<\varepsilon$.
\end{enumerate}
Then $\mathcal{Y}+\mathcal{Z}$ is closed in $\mathcal{X}$.
\end{theorem}

Let $1<p<\infty$ and $M_p$ be the Banach algebra of symbols $a$ generating
bounded Laurent operators $L(a)$ on $\ell^p(\mathbb{Z})$, let $C_p$ be the
closure of trigonometric polynomials in $M_p$ and 
$H_p^\infty=H^\infty\cap M_p$. In 1987, B\"ottcher and Silbermann proved 
that $C_p+H_p^\infty$ is a closed subalgebra of $M_p$ 
(see \cite[Theorem~1]{BS89} and also \cite[Theorem~2.53]{BS06}).
Their proof is based on Theorem~\ref{th:Zalcman-Rudin} with
$\mathcal{X}=M_p$, $\mathcal{Y}=C_p$, $\mathcal{Z}=H_p^\infty$, and
$\Phi$ being the collection of all convolution operators with the
Fej\'er kernels.

The aim of this paper is to extend the 
{result by B\"ottcher and Silbermann mentioned above}
to the setting of reflexive 
rearrangement-invariant Banach sequence spaces (see \cite[Ch.~2]{BS88}
and Section~\ref{subsec:riBSS} {below}) 
with nontrivial Boyd indices $\alpha_X,\beta_X$
(see \cite{B69} and Section~\ref{sec:Boyd}) 
and symmetric weights belonging to the intersection of the Muckenhoupt 
classes $A_{1/\alpha_X}(\mathbb{Z})$ and $A_{1/\beta_X}(\mathbb{Z})$. 
Note that the class of 
rearrangement-invariant Banach sequence spaces is very wide. The most
interesting examples are Lebesgue sequence spaces $\ell^p(\mathbb{Z})$
with $1\le p\le\infty$, Lorentz sequence spaces $\ell^{p,q}(\mathbb{Z})$
with $1<p<\infty$ and $1\le q\le\infty$, Orlicz sequence spaces 
$\ell^\Phi(\mathbb{Z})$ with suitable convex functions 
$\Phi:[0,\infty)\to[0,\infty]$. Boyd indices $\alpha_X$ and $\beta_X$
satisfy $0\le\alpha_X\le\beta_X\le 1$ and, roughly speaking, play a role
of $1/p$ for the spaces $\ell^p(\mathbb{Z})$. More precisely,
$\alpha_{\ell^p}=\beta_{\ell^p}=1/p$ for all $1\le p\le\infty$,
however there are reflexive Orlicz spaces $\ell^\Phi(\mathbb{Z})$
such that $0<\alpha_{\ell^\Phi}<\beta_{\ell^\Phi}<1$
(see \cite[Corollary~11.12, Theorem~11.5]{M89} and \cite{B71}).
It is customary to say that Boyd indices of a rearrangement-invariant 
space $X(\mathbb{Z})$ are nontrivial if
\[
0<\alpha_X,\quad \beta_X<1.
\]
Without going into technical details, we mention that our main result 
(Theorem~\ref{th:algebra-C-plus-H-infinity}) {seems to be} new already 
in the setting of weighted Lebesgue sequence spaces $\ell^p(\mathbb{Z},w)$ 
with $1<p<\infty$ and symmetric Muckenhoupt weights $w\in A_p(\mathbb{Z})$.

Now we give some important definitions.
Let $X(\mathbb{Z})$ be a Banach sequence space (non necessarily 
rearrangement-invariant) (see \cite[Ch.~1]{BS88} and Section~\ref{subsec:BSS}).
Let $\mathcal{P}'$ be the space of periodic distributions (see, e.g.,
\cite[Ch.~3 and~5]{B73} and Section~\ref{sec:distributions})
and let 
$S_0(\mathbb{Z})$ denote the set of all finitely supported sequences.
For $a\in\mathcal{P}'$ and $\varphi\in S_0(\mathbb{Z})$, we define
the convolution $a*\varphi$ as the sequence
\[
(a*\varphi)_j=\sum_{k\in\mathbb{Z}}\widehat{a}(j-k)\varphi_k,
\quad
j\in\mathbb{Z},
\]
where $\{\widehat{a}(j)\}_{j\in\mathbb{Z}}$ is the sequence of Fourier
coefficients of the distribution $a$. By $M_{X(\mathbb{Z})}$ we denote
the collection of all distributions $a\in\mathcal{P}'$ for which 
$a*\varphi\in X(\mathbb{Z})$ whenever $\varphi\in S_0(\mathbb{Z})$ and
\begin{equation}\label{eq:multiplier-norm}
\|a\|_{M_{X(\mathbb{Z})}}:=\sup\left\{
\frac{\|a*\varphi\|_{X(\mathbb{Z})}}{\|\varphi\|_{X(\mathbb{Z})}}\ :\
\varphi\in S_0(\mathbb{Z}), \ \varphi\ne 0
\right\}<\infty.
\end{equation}

If the space $X(\mathbb{Z})$ is separable, then $S_0(\mathbb{Z})$ is dense 
in $X(\mathbb{Z})$ (see Lemma~\ref{le:density}(a) below). 
In this case, for $a\in M_{X(\mathbb{Z})}$, the operator
from $S_0(\mathbb{Z})$ to $X(\mathbb{Z})$ defined by $\varphi\mapsto a*\varphi$
extends to a bounded operator 
\[
L(a):X(\mathbb{Z})\to X(\mathbb{Z}),
\quad
\varphi\mapsto a*\varphi,
\]
which is referred to as the Laurent operator with symbol $a$.

Let $\mathbb{Z}_+:=\{0,1,2,\dots\}=\{0\}\cup\mathbb{N}$. We will equip 
$\mathbb{S}\in\{\mathbb{N},\mathbb{Z}_+,\mathbb{Z}\}$ with the 
counting measure $m$, that is, $m$ is a purely atomic measure and $m(\{k\})=1$ 
for every $k\in\mathbb{S}$.

A weight on $\mathbb{S}$ is a sequence $w=\{w_k\}_{k\in\mathbb{S}}$ of 
positive numbers. Let $l,n\in\mathbb{S}$ satisfy $l\le n$. We call the set 
of the form $J=\{l,\dots,n\}$ an interval of $\mathbb{S}$. Let $1<p<\infty$ 
and $1/p+1/q=1$. A weight $w=\{w_k\}_{k\in\mathbb{S}}$ is said to belong to 
the Muckenhoupt class $A_p(\mathbb{S})$ if 
\[
[w]_{A_p(\mathbb{S})}:= \sup_{J\subset \mathbb{S}} 
\frac{1}{m(J)}
\left(\sum_{k\in J} w_k^p\right)^{1/p}
\left(\sum_{k\in J} w_k^{-q}\right)^{1/q}  
< \infty,
\]
where the supremum is taken over all intervals $J\subset\mathbb{S}$.

A weight $w=\{w_k\}_{k\in\mathbb{Z}}$ is said to be symmetric if $w_{-k}=w_k$ 
for all $k\in\mathbb{Z}$. The collection of all symmetric weights in the 
Muckenhoupt class $A_p(\mathbb{Z})$ will be denoted by 
$A_p^{\rm sym}(\mathbb{Z})$.

For a Banach sequence
space $X(\mathbb{Z})$ and a weight $w=\{w_k\}_{k\in\mathbb{Z}}$,
the weighted Banach sequence space $X(\mathbb{Z},w)$ consists of 
all sequences $f=\{f_k\}_{k\in\mathbb{Z}}$ such that
$fw:=\{f_kw_k\}_{k\in\mathbb{Z}}$ belongs to $X(\mathbb{Z})$. It is easy to 
see that $X(\mathbb{Z},w)$ is itself a Banach sequence space with respect
to the norm
\[
\|f\|_{X(\mathbb{Z},w)}:=\|fw\|_{X(\mathbb{Z})}.
\]

For $1\le p\le\infty$, let $L^p(-\pi,\pi)$ denote the space of
all $2\pi$-periodic measurable functions $f:\mathbb{R}\to\mathbb{C}$ 
such that
\[
\|f\|_{L^p(-\pi,\pi)}
:=
\left\{\begin{array}{ll}
\displaystyle
\left(\frac{1}{2\pi}\int_{-\pi}^\pi |f(\theta)|^p\,d\theta\right)^{1/p},
&
\quad 1\le p<\infty,
\\[3mm]
\displaystyle
\operatornamewithlimits{ess\,sup}_{\theta\in(-\pi,\pi)}|f(\theta)|,
&
\quad p=\infty,
\end{array}\right.
\]
is finite. A $2\pi$-periodic function $a:\mathbb{R}\to\mathbb{C}$ is said to 
be of bounded variation if
\[
V(a):=\sup\sum_{j=1}^n |a(\theta_{j+1})-a(\theta_j)|<\infty,
\]
where the supremum is taken over all partitions 
$-\pi\le\theta_1<\theta_2<\dots<\theta_{n+1}=\pi$. It is well known
that the collection of all functions of bounded variation $BV[-\pi,\pi]$
is a Banach algebra with respect to the norm
\[
\|a\|_{BV[-\pi,\pi]}:=\|a\|_{L^\infty(-\pi,\pi)}+V(a).
\]

Our first result is a generalization of properties of $M_{X(\mathbb{Z},w)}$
listed for $\ell^p(\mathbb{Z})$ and $w=1$ in \cite[Section~2.5(d), (f), (g)]{BS06}.
\begin{theorem}\label{th:Banach-algebra-of-multipliers}
Let $X(\mathbb{Z})$ be a reflexive rearrangement-invariant Banach sequence
space with nontrivial Boyd indices $\alpha_X,\beta_X$ and let
$w\in A_{1/\alpha_X}^{\rm sym}(\mathbb{Z})\cap 
A_{1/\beta_X}^{\rm sym}(\mathbb{Z})$.  
\begin{enumerate}
\item[(a)]
If $a\in M_{X(\mathbb{Z},w)}$, then $a\in L^\infty(-\pi,\pi)$ and
\[
\|a\|_{L^\infty(-\pi,\pi)}\le \|a\|_{M_{X(\mathbb{Z},w)}}.
\]

\item[(b)]
$M_{X(\mathbb{Z},w)}$ is a Banach algebra under pointwise multiplication 
and the norm
\[
\|a\|_{M_{X(\mathbb{Z},w)}}=\|L(a)\|_{\mathcal{B}(X(\mathbb{Z},w))}.
\]

\item[(c)] 
There is a constant $C_{X(\mathbb{Z},w)}\in(0,\infty)$ depending 
only on $X(\mathbb{Z})$ and $w$ such that
\[
\|a\|_{M_{X(\mathbb{Z},w)}}
\le 
C_{X(\mathbb{Z},w)}\left(\|a\|_{L^\infty(-\pi,\pi)}+V(a)\right)
\] 
for all $a\in BV[-\pi,\pi]$. In particular, 
$BV[-\pi,\pi]\hookrightarrow M_{X(\mathbb{Z},w)}$.
\end{enumerate}
\end{theorem}

It would be interesting to clarify whether the above result remains true
if one drops the assumption that the weight $w$ is symmetric. 

Let $\mathcal{TP}$ denote the set of all trigonometric polynomials,
that is, functions of the form
\[
p(\theta)=\sum_{|k|\le n} c_k e^{ik\theta},
\quad\theta\in\mathbb{R},
\]
where $n\in\mathbb{Z}_+$ and $c_k\in\mathbb{C}$ for all $k\in\{-n,\dots,n\}$.
Since each function in $\mathcal{TP}$ is of bounded variation, it follows 
from Theorem~\ref{th:Banach-algebra-of-multipliers}(c)
that $\mathcal{TP}\subset M_{X(\mathbb{Z},w)}$. We denote by 
$C_{X(\mathbb{Z},w)}$ the closure of $\mathcal{TP}$ with respect 
to the norm of the Banach algebra $M_{X(\mathbb{Z},w)}$ and put
\[
H_{X(\mathbb{Z},w)}^{\infty,\pm}:=
\left\{a\in M_{X(\mathbb{Z},w)}:\widehat{a}(\pm n)=0 \mbox{ for }n<0\right\}.
\]

The following extension of \cite[Theorem~1]{BS89} and \cite[Theorem~2.53]{BS06} 
is the main result of this paper.
\begin{theorem}
\label{th:algebra-C-plus-H-infinity}
Let $X(\mathbb{Z})$ be a reflexive rearrangement-invariant Banach sequence
space with nontrivial Boyd indices $\alpha_X,\beta_X$ and let
$w\in A_{1/\alpha_X}^{\rm sym}(\mathbb{Z})\cap 
A_{1/\beta_X}^{\rm sym}(\mathbb{Z})$. Then the sets 
$C_{X(\mathbb{Z},w)}+H_{X(\mathbb{Z},w)}^{\infty,\pm}$ 
are closed subalgebras of the Banach algebra $M_{X(\mathbb{Z},w)}$.
\end{theorem}
The paper is organized as follows. The first two sections 
{cover the} preliminaries. 
In Section~\ref{sec:BSS}, we collect 
properties of Banach sequence spaces culminating with the fact that
if an operator is bounded on a Banach sequence space $X(\mathbb{Z})$
and on its associate space $X'(\mathbb{Z})$, then it is bounded on the 
space $\ell^2(\mathbb{Z})$. Section~\ref{sec:riBSS} is devoted to 
rearrangement-invariant Banach sequence spaces, their Boyd indices, and
the Boyd interpolation theorem \cite{B69}. 

In Section~\ref{sec:Muckenhoupt-stability}, we prove that if 
$w\in A_{p_0}^{\rm sym}(\mathbb{Z})$, then 
$w^{1+\varepsilon}\in A_p^{\rm sym}(\mathbb{Z})$ for all $\varepsilon$
sufficiently close to zero and all $p$ sufficiently close to $p_0$. This 
result extends the previous stability result for discrete Muckenhoupt 
weights by B\"ottcher and Seybold \cite[Theorem~1.1]{BSey99}, which did 
not allow to ``shake" $p_0$. We follow closely the proof
given by B\"ottcher and Karlovich \cite[Theorem~2.31]{BK97} for continuous
Muckenhoupt weights. 

Section~\ref{sec:algebra-MX} is dedicated to the study of properties of
$M_{X(\mathbb{Z})}$. We show that if $X(\mathbb{Z})$ is a 
reflection-invariant Banach sequence space, then 
$M_{X(\mathbb{Z})}=M_{X'(\mathbb{Z})}$ with equal norms. If, in addition, 
$X(\mathbb{Z})$ is reflexive, then the interpolation theorem from
Section~\ref{sec:BSS} allows us to prove that 
$M_{X(\mathbb{Z})}\hookrightarrow L^\infty(-\pi,\pi)$ and that
$M_{X(\mathbb{Z})}$ is a Banach algebra. We conclude this section
with the proof of the uniform boundedness of Fej\'er means of 
$a\in M_{X(\mathbb{Z})}$ if $X(\mathbb{Z})$ is a reflexive 
reflection-invariant Banach sequence space. 

We start Section~\ref{sec:proofs} with the observation that if 
$X(\mathbb{Z})$ is a reflexive rearrangement-invariant Banach sequence 
space and a weight $w$ is symmetric, then $X(\mathbb{Z},w)$ is a 
reflexive reflection-invariant Banach sequence space. So, 
Theorem~\ref{th:Banach-algebra-of-multipliers}(a)--(b) is a consequence of
results of Section~\ref{sec:algebra-MX}. Further, combining the stability
result for Muckenhoupt weights of Section~\ref{sec:Muckenhoupt-stability},
the Boyd interpolation theorem (Theorem~\ref{th:Boyd}), and Stechkin's
inequality for $\ell^p(\mathbb{Z},w)$ with $w\in A_p(\mathbb{Z})$
proved by B\"ottcher and Seybold \cite[Theorem~1.2]{BSey99}, we prove
Theorem~\ref{th:Banach-algebra-of-multipliers}(c).

The proof of Theorem~\ref{th:algebra-C-plus-H-infinity} essentially 
consists {of} verification of properties (i)--(iv) of 
Theorem~\ref{th:Zalcman-Rudin} for 
$\mathcal{X}=M_{X(\mathbb{Z},w)}$,
$\mathcal{Y}=C_{X(\mathbb{Z},w)}$, 
$\mathcal{Z}:=H_{X(\mathbb{Z},w)}^{\infty,\pm}$
and $\Phi$ being the collection of all convolution operators with the
Fej\'er kernels $a\mapsto\sigma_n(a)$. The most difficult is the 
verification of property (iv), where we have to prove that $\sigma_n(a)\to a$
as $n\to\infty$  in the norm of $M_{X(\mathbb{Z},w)}$ for every
$a\in C_{X(\mathbb{Z},w)}$. It follows from the definition of 
$C_{X(\mathbb{Z},w)}$ and the last result of Section~\ref{sec:algebra-MX}
that it is enough to consider 
$a\in\mathcal{TP}\subset C[-\pi,\pi]\cap BV[-\pi,\pi]$. On the other hand, 
the combination of the Boyd and Stein-Weiss interpolation theorems with
the Stechkin inequality allows us to to reduce the study of convergence 
$\sigma_n(a)\to a$ in the norm of $M_{X(\mathbb{Z},w)}$ for 
$a\in C[-\pi,\pi]\cap BV[-\pi,\pi]$ to the study
of convergence $\sigma_n(a)\to a$ in the norm of 
$M_{\ell^2(\mathbb{Z})}=L^\infty(-\pi,\pi)$, which is a classical result.
Once again, our ability to perform this interpolation argument is based on the
stability result for symmetric Muckenhoupt weights obtained in 
Section~\ref{sec:Muckenhoupt-stability}.
This will complete the proof of Theorem~\ref{th:algebra-C-plus-H-infinity}.
\section{Preliminaries on Banach sequence spaces}\label{sec:BSS}
\subsection{Banach sequence spaces}\label{subsec:BSS}
Let $\mathbb{S}\in\{\mathbb{N},\mathbb{Z}\}$, let $\ell^0(\mathbb{S})$ be the 
linear space of all sequences $f:\mathbb{S}\to\mathbb{C}$, and let 
$\ell_+^0(\mathbb{S})$ be the cone of nonnegative sequences in 
$\ell^0(\mathbb{S})$. According to 
\cite[Ch.~1, Definition~1.1]{BS88}, a Banach function norm 
$\varrho:\ell_+^0(\mathbb{S})\to [0,\infty]$ is a mapping which satisfies 
the following axioms for all $f,g\in \ell_+^0(\mathbb{S})$, for all 
{sequences}
$\{f^{(n)}\}_{n\in\mathbb{N}}$ in $\ell_+^0(\mathbb{S})$, for all finite 
subsets $E\subset\mathbb{S}$, and all constants $\alpha\ge 0$:
\begin{eqnarray*}
{\rm (A1)} & &
\varrho(f)=0  \Leftrightarrow  f=0,
\
\varrho(\alpha f)=\alpha\varrho(f),
\
\varrho(f+g) \le \varrho(f)+\varrho(g),\\
{\rm (A2)} & &0\le g \le f \ \ \Rightarrow \ 
\varrho(g) \le \varrho(f)
\quad\mbox{(the lattice property)},\\
{\rm (A3)} & &0\le f^{(n)} \uparrow f \ \ \Rightarrow \
       \varrho(f^{(n)}) \uparrow \varrho(f)\quad\mbox{(the Fatou property)},\\
{\rm (A4)} & & \varrho(\chi_E) <\infty,\\
{\rm (A5)} & &\sum_{k\in E} f_k \le C_E\varrho(f) ,
\end{eqnarray*}
where $\chi_E$ is the characteristic (indicator) function of $E$, and the 
constant $C_E \in (0,\infty)$ may depend on $\varrho$ and $E$, 
but is independent of $f$. The set $X(\mathbb{S})$ of all sequences 
$f\in \ell^0(\mathbb{S})$ for which $\varrho(|f|)<\infty$ is called a 
Banach sequence space. For each $f\in X(\mathbb{S})$, the norm of $f$ is 
defined by 
\[
\|f\|_{X(\mathbb{S})} :=\varrho(|f|). 
\]
The set $X(\mathbb{S})$ equipped with the natural linear space operations 
and this norm becomes a Banach space 
(see \cite[Ch.~1, Theorems~1.4 and~1.6]{BS88}). If $\varrho$ 
is a Banach function norm, its associate norm $\varrho'$ is defined on 
$\ell_+^0(\mathbb{S})$ by
\[
\varrho'(g):=\sup\left\{
\sum_{k\in\mathbb{S}} f_kg_k \ : \ 
f=\{f_k\}_{k\in\mathbb{S}}\in\ell_+^0(\mathbb{S}), \ \varrho(f) \le 1
\right\}, \quad g\in \ell_+^0(\mathbb{S}).
\]
It is a Banach function norm itself \cite[Ch.~1, Theorem~2.2]{BS88}.
The Banach sequence space $X'(\mathbb{S})$ determined by the Banach 
function norm $\varrho'$ is called the associate space (K\"othe dual) 
of $X(\mathbb{S})$. The associate space $X'(\mathbb{S})$ can be viewed 
as a subspace of the Banach dual space $X^*(\mathbb{S})$.

Let us conclude this subsections with two results highlighting the idea
that in many {cases} arbitrary sequences can be approximated 
by finitely supported sequences.
\begin{lemma}\label{le:density}
Let $X(\mathbb{Z})$ be a Banach sequence space and $X'(\mathbb{Z})$ be
its associate space.
\begin{enumerate}
\item[(a)]
If $X(\mathbb{Z})$ is separable, then the set $S_0(\mathbb{Z})$
is dense in the space $X(\mathbb{Z})$.

\item[(b)]
If $X(\mathbb{Z})$ is reflexive, then the set $S_0(\mathbb{Z})$ is
dense in the space $X(\mathbb{Z})$ and in the space $X'(\mathbb{Z})$.
\end{enumerate}
\end{lemma}
Part (a) follows from 
\cite[Ch.~1, Proposition~3.10, Theorem 3.11, and Corollary~5.6]{BS88}.
Part (b) follows from part (a) and \cite[Ch.~1, Corollary~4.4]{BS88}.
\begin{lemma}[{see \cite[Lemma~2.1]{KS24-PAFA} and also 
\cite[Lemma~2.10]{KS19}}]
\label{le:norm-on-nice-sequences}
Let $X(\mathbb{Z})$ be a Banach sequence space and $X'(\mathbb{Z})$ 
be its associate space. For every sequence
$f=\{f_k\}_{k\in\mathbb{Z}}\in X(\mathbb{Z})$,
\[
\|f\|_{X(\mathbb{Z})}
=
\sup\left\{\left|\sum_{k\in\mathbb{Z}} f_ks_k\right|\ :\
s=\{s_k\}_{k\in\mathbb{Z}}\in S_0(\mathbb{Z}),\ 
\|s\|_{X'(\mathbb{Z})}\le 1\right\}.
\]
\end{lemma}
\subsection{Minkowski integral inequality for Banach sequence spaces}
We will need the following Minkowski integral inequality for Banach 
sequence spaces, which is a special case of \cite[Proposition~2.1]{S95}.
\begin{lemma}\label{le:Minkowski-inequality}
Let $\mathbb{Z}$ be equipped with the counting measure $m$ and $[-\pi,\pi]$
be equipped with a $\sigma$-finite probability measure $\mu$. 
Suppose $X(\mathbb{Z})$ is a Banach sequence space.
If $f$ is $m\times\mu$-measurable, then
\[
\left\|
\left\{\int_{-\pi}^\pi |f(j,x)|\,d\mu(x)\right\}_{j\in\mathbb{Z}}
\right\|_{X(\mathbb{Z})}
\le 
\int_{-\pi}^\pi 
\left\|\{f(j,x)\}_{j\in\mathbb{Z}}
\right\|_{X(\mathbb{Z})}d\mu(x).
\]
\end{lemma}
\subsection{Reflection-invariant Banach sequence spaces}
By analogy with \cite[Subsection~7.2]{KS19}, we say that a Banach sequence 
space $X(\mathbb{Z})$ is reflection-invariant if 
$\|\varphi\|_{X(\mathbb{Z})}=\|\widetilde{\varphi}\|_{X(\mathbb{Z})}$
for every $\varphi\in X(\mathbb{Z})$, where $\widetilde{\varphi}$
denotes the reflection of a sequence $\varphi=\{\varphi_k\}_{k\in\mathbb{Z}}$
defined by
\[
\widetilde{\varphi}_k:=\varphi_{-k},
\quad
k\in\mathbb{Z}.
\]
\begin{lemma}\label{le:reflection-invariant-spaces}
A Banach sequence space $X(\mathbb{Z})$ is reflection-invariant if and only 
if its associate space $X'(\mathbb{Z})$ is reflection-invariant. 
\end{lemma}
\begin{proof}
Fix $\psi=\{\psi_k\}_{k\in\mathbb{Z}}\in X'(\mathbb{Z})$. If $X(\mathbb{Z})$
is reflection-invariant, then it follows from \cite[Ch.~1, Lemma~2.8]{BS88}
that 
\begin{align*}
\|\widetilde{\psi}\|_{X'(\mathbb{Z})}
&=
\sup\left\{\left|\sum_{k\in\mathbb{Z}}\widetilde{\psi}_k\varphi_k\right|
\ :\ \varphi=\{\varphi_k\}_{k\in\mathbb{Z}}\in X(\mathbb{Z}),\
\|\varphi\|_{X(\mathbb{Z})}\le 1\right\}
\\
&=
\sup\left\{\left|\sum_{k\in\mathbb{Z}}\psi_k\widetilde{\varphi}_k\right|
\ :\ \varphi=\{\varphi_k\}_{k\in\mathbb{Z}}\in X(\mathbb{Z}),\
\|\varphi\|_{X(\mathbb{Z})}\le 1\right\}
\\
&=
\sup\left\{\left|\sum_{k\in\mathbb{Z}}\psi_k\varphi_k\right|
\ :\ \varphi=\{\varphi_k\}_{k\in\mathbb{Z}}\in X(\mathbb{Z}),\
\|\varphi\|_{X(\mathbb{Z})}\le 1\right\}
=
\|\psi\|_{X'(\mathbb{Z})}
\end{align*}
because $\|\widetilde{\varphi}\|_{X(\mathbb{Z})}=\|\varphi\|_{X(\mathbb{Z})}$
for all $\varphi=\{\varphi_k\}_{k\in\mathbb{Z}}\in X(\mathbb{Z})$.
Hence $X'(\mathbb{Z})$ is reflection-invariant.

If $X'(\mathbb{Z})$ is reflection-invariant, then, by what has already been 
proved, $X''(\mathbb{Z})$ is reflection-invariant. It remains to recall that,
by the Lorentz-Luxemburg theorem (see \cite[Ch.~1, Theorem~2.7]{BS88}),
the space 
{
$X''(\mathbb{Z})$ coincides with the space $X(\mathbb{Z})$ and
$\|\varphi\|_{X(\mathbb{Z})}=\|\varphi\|_{X''(\mathbb{Z})}$ 
for all  $\varphi=\{\varphi_k\}_{k\in\mathbb{Z}}\in X(\mathbb{Z})$.
}
Thus, $X(\mathbb{Z})$ is reflection-in\-vari\-ant.
\end{proof}
\subsection{Interpolation of operators between a Banach sequence 
space and its associate space}
Let $X_0(\mathbb{Z})$ and $X_1(\mathbb{Z})$ be Banach sequence spaces. Suppose
$0<\theta<1$. The Calder\'on product $(X_0^{1-\theta} X_1^\theta)(\mathbb{Z})$
of $X_0(\mathbb{Z})$ and $X_1(\mathbb{Z})$ consists of all sequences
$x=\{x_k\}_{k\in\mathbb{Z}}\in\ell^0(\mathbb{Z})$ such that the inequality
\begin{equation}\label{eq:Calderon-product}
|x_k|\le\lambda |y_k|^{1-\theta}|z_k|^\theta,
\quad
k\in\mathbb{Z},
\end{equation}
holds for some $\lambda>0$ and elements 
$y=\{y_k\}_{k\in\mathbb{Z}}\in X_0(\mathbb{Z})$ and
$z=\{z_k\}_{k\in\mathbb{Z}}\in X_1(\mathbb{Z})$ with
$\|y\|_{X_0(\mathbb{Z})}\le 1$ and $\|z\|_{X_1(\mathbb{Z})}\le 1$
(see \cite{C64}).
The norm of $x$ in $(X_0^{1-\theta} X_1^\theta)(\mathbb{Z})$
is defined to be the infimum of all values of $\lambda$ in
\eqref{eq:Calderon-product}.
\begin{theorem}\label{th:Calderon-Lozanovskii}
Let $X(\mathbb{Z})$ be a Banach sequence space and $X'(\mathbb{Z})$ be its
associate space. If $A:\ell^0(\mathbb{Z})\to\ell^0(\mathbb{Z})$ 
{is} a linear
operator bounded on $X(\mathbb{Z})$ and on $X'(\mathbb{Z})$, then $A$
is bounded on $\ell^2(\mathbb{Z})$ and
\[
\|A\|_{\mathcal{B}(\ell^2(\mathbb{Z}))}
\le
\|A\|_{\mathcal{B}(X(\mathbb{Z}))}^{1/2}
\|A\|_{\mathcal{B}(X'(\mathbb{Z}))}^{1/2}.
\]
\end{theorem} 
\begin{proof}
It follows from the Lozanovskii interpolation formula that
\begin{equation}\label{eq:Lozanovskii}
(X^{1/2}(X')^{1/2})(\mathbb{Z})=\ell^2(\mathbb{Z})
\end{equation}
with equality of the norms (see, e.g., \cite[Ch.~15, Example~6, p.~185]{M89}).
The desired result is a consequence of \eqref{eq:Lozanovskii} and
an abstract version of the Riesz-Thorin interpolation theorem given in
\cite[Theorem~3.11]{MS19}.
\end{proof}
\section{Rearrangement-invariant Banach sequence spaces\\ and their Boyd indices}
\label{sec:riBSS}
\subsection{Rearrangement-invariant Banach sequence spaces}\label{subsec:riBSS}
Let $\mathbb{S}\in\{\mathbb{N},\mathbb{Z}\}$. The distribution function of a 
sequence $f=\{f_k\}_{k\in\mathbb{S}}\in\ell^0(\mathbb{S})$ is defined by
\[
d_f(\lambda):=m\{k\in\mathbb{S}:|f_k|>\lambda\},
\quad
\lambda\ge 0,
\]
where $m(S)$ denotes the measure (cardinality) of a set $S\subset\mathbb{S}$.
For any sequence $f\in\ell^0(\mathbb{S})$, its decreasing rearrangement
is given by
\[
f^*(n):=\inf\{\lambda\ge 0:d_f(\lambda)\le n-1\},
\quad
n\in\mathbb{N}.
\]
One says that sequences 
$f=\{f_k\}_{k\in\mathbb{S}_1}\in\ell^0(\mathbb{S}_1)$, 
$g=\{g_k\}_{k\in\mathbb{S}_2}\in\ell^0(\mathbb{S}_2)$ 
with $\mathbb{S}_1,\mathbb{S}_2\in\{\mathbb{N},\mathbb{Z}\}$
are equimeasurable if $d_f=d_g$. 

A Banach function norm $\varrho:\ell_+^0(\mathbb{S})\to[0,\infty]$ is said to 
be rearrangement-invariant if $\varrho(f)=\varrho(g)$ for every pair of 
equimeasurable sequences 
$f=\{f_k\}_{k\in\mathbb{S}}$, $g=\{g_k\}_{k\in\mathbb{S}}\in\ell_+^0(\mathbb{S})$. 
In that case, the Banach sequence space $X(\mathbb{S})$ generated by $\varrho$ 
is said to be a rearrangement-invariant Banach sequence space 
(cf. \cite[Ch.~2, Definition~4.1]{BS88}).
It follows from \cite[Ch.~2, Proposition~4.2]{BS88} that if a Banach sequence 
space $X(\mathbb{S})$ is rearrangement-invariant, then its associate space 
$X'(\mathbb{S})$ is also a rearrangement-invariant Banach sequence space.
\subsection{Boyd indices}\label{sec:Boyd}
Let $X(\mathbb{Z})$ be a rearrangement-invariant Banach sequence
space generated by a rearrangement-invariant Banach function norm
$\varrho:\ell_+^0(\mathbb{Z})\to[0,\infty]$. By the Luxemburg 
representation theorem (see \cite[Ch.~2, Theorem~4.10]{BS88}), there 
exists a unique rearrangement-invariant Banach function norm 
$\overline{\varrho}:\ell_+^0(\mathbb{N})\to[0,\infty]$
such that
\[
\varrho(f)=\overline{\varrho}(f^*),
\quad f\in\ell_+^0(\mathbb{Z}).
\]
The rearrangement-invariant Banach sequence space generated by 
$\overline{\varrho}$ is denoted by $\overline{X}(\mathbb{N})$ and is 
called the Luxemburg representation of $X(\mathbb{Z})$.

Let $t\ge 0$. We denote by $\lfloor t\rfloor$ the greatest integer less than 
or equal to $t$. For $f=\{f_k\}_{k\in\mathbb{N}}:\mathbb{N}\to[0,\infty)$, 
consider the following operators:
\[
(E_jf)_k:=f_{jk},
\quad
(F_jf)_k:=f_{\lfloor (k-1)/j\rfloor+1},
\quad
j,k\in\mathbb{N}.
\]
For $j\in\mathbb{N}$, let
\begin{align*}
H(j,X)
&:=
\sup\left\{\overline{\varrho}(E_j f^*)\ :\ f\in \overline{X}(\mathbb{N}),\
\overline{\varrho}(f)\le 1\right\},
\\
K(j,X)
&:=
\sup\left\{\overline{\varrho}(F_j f^*)\ :\ f\in \overline{X}(\mathbb{N}),\
\overline{\varrho}(f)\le 1\right\}.
\end{align*}
It follows from \cite[Lemmas~3--5]{B69} that the following limits
\[
\alpha_X:=\lim_{j\to\infty}\frac{-\log H(j,X)}{\log j},
\quad
\beta_X:=\lim_{j\to\infty}\frac{\log K(j,X)}{\log j}
\]
exist and satisfy
\[
0\le\alpha_X\le\beta_X\le 1,
\quad
\alpha_{X'}=1-\beta_X,
\quad
\beta_{X'}=1-\alpha_X.
\]
The numbers $\alpha_X$ and $\beta_X$ are called the lower and the upper 
Boyd indices, respectively. 
\subsection{Boyd's interpolation theorem}
The next theorem follows from the Boyd interpolation theorem 
\cite[Theorem~1]{B69} for quasi-linear operators of weak types $(p,p)$ and
$(q,q)$. Its proof for nonatomic measure spaces can also be found in 
\cite[Ch.~3, Theorem~5.16]{BS88} and \cite[Theorem~2.b.11]{LT79}.
\begin{theorem}\label{th:Boyd}
Let $1\le q<p\le\infty$ and $X(\mathbb{Z})$ be a rearrangement-invariant 
Banach sequence space with the Boyd indices $\alpha_X,\beta_X$ satisfying
\[
1/p<\alpha_X\le \beta_X<1/q. 
\]
Then there exists a constant $C_{p,q}\in(0,\infty)$ with the following 
property. If a linear operator $T:\ell^0(\mathbb{Z})\to\ell^0(\mathbb{Z})$
is bounded on the spaces $\ell^p(\mathbb{Z})$ and $\ell^q(\mathbb{Z})$, 
then it is  also bounded on the rearrangement-invariant Banach sequence 
space $X(\mathbb{Z})$ and
\begin{equation}\label{eq:Boyd}
\|T\|_{\mathcal{B}(X(\mathbb{Z}))}\le C_{p,q}\max
\big\{
\|T\|_{\mathcal{B}(\ell^p(\mathbb{Z}))},
\|T\|_{\mathcal{B}(\ell^q(\mathbb{Z}))}
\big\}.
\end{equation}
\end{theorem}
Notice that estimate \eqref{eq:Boyd} is not stated explicitly in
\cite{BS88,B69,LT79}. However, it can be extracted from the proof of the
Boyd interpolation theorem.
\section{Discrete Muckenhoupt weights}\label{sec:Muckenhoupt-stability}
\subsection{Symmetric reproduction of Muckenhoupt weights}
{
Recall that $A_p^{\rm sym}(\mathbb{Z})$ denotes the call of weights
$v\in A_p(\mathbb{Z})$ satisfying $v_{-k}=v_k$ for all $k\in\mathbb{Z}$.}
The following lemma allows us to extend a weight $w\in A_p(\mathbb{Z}_+)$
to a weight $v\in A_p^{\rm sym}(\mathbb{Z})$ (see \cite[Lemma~3.2]{KS24-LAA}
and also \cite[Section~2.4]{BK97}, where similar
results in the continuous case are considered).
\begin{lemma}\label{le:even-extension} 
Let $1 < p < \infty$, $w : \mathbb{Z}_+ \to (0, \infty)$ and $v_n := w_{|n|}$, 
$n \in \mathbb{Z}$. Then 
\[
w \in A_p(\mathbb{Z}_+) 
\quad \Longleftrightarrow \quad 
v \in A_p^{\rm sym}(\mathbb{Z}).
\]
\end{lemma}
\subsection{Reverse H\"older inequality}
We will need the following result contained in \cite[Lemma~2.3]{BSey99}.
\begin{theorem}\label{th:RHI}
Let $1<p<\infty$ and $w\in A_p(\mathbb{Z}_+)$. Then there exist $\delta>0$
and $C>0$ depending only on $p$ and $w$ such that
\[
\left(\frac{1}{m(R)}\sum_{k\in R}w_k^{p(1+\delta)}\right)^{1/(1+\delta)}
\le 
\frac{C}{m(R)}\sum_{k\in R}w_k^p
\]
for all intervals $R\subset\mathbb{Z}_+$ {such that} 
$m(R)=2^r$ with $r\in\mathbb{N}$.
\end{theorem}
\subsection{``Convexity" property of discrete Muckenhoupt weights}
The following lemma is analogous to \cite[Proposition~2.30]{BK97}. We give
its proof for the sake of completeness of presentation.
\begin{lemma}\label{le:BK230}
Let $w:\mathbb{Z}_+\to(0,\infty)$ be a weight. Then the set
\[
\Gamma:=\left\{
(p,\delta)\in(1,\infty)\times\mathbb{R}\ :\ w^{\delta/p}\in A_p(\mathbb{Z}_+)
\right\} 
\]
is convex.
\end{lemma}
\begin{proof}
Let $(p_1,\delta_1)\in\Gamma$ and $(p_2,\delta_2)\in\Gamma$ and put
\[
p(\theta):=(1-\theta)p_1+\theta p_2,
\quad
\delta(\theta):=(1-\theta)\delta_1+\theta\delta_2,
\quad
\theta\in(0,1).
\]
By definition of $\Gamma$, we have $w^{\delta_1/p_1}\in A_{p_1}(\mathbb{Z}_+)$
and $w^{\delta_2/p_2}\in A_{p_2}(\mathbb{Z}_+)$. We have to show that
$w^{\delta(\theta)/p(\theta)}\in A_{p(\theta)}(\mathbb{Z}_+)$ for all
$\theta\in(0,1)$. 

Let $J\subset\mathbb{Z}_+$ be any interval. Consider
\begin{equation}\label{eq:BK230-1}
\alpha:=1/(1-\theta),
\quad
\beta:=1/\theta.
\end{equation}
Then for all $k\in J$,
\[
\left(w_k^{\delta(\theta)/p(\theta)}\right)^{p(\theta)}
=
w_k^{(1-\theta)\delta_1+\theta\delta_2}
=
\left(w_k^{\delta_1}\right)^{1/\alpha}
\left(w_k^{\delta_2}\right)^{1/\beta}.
\]
H\"older's inequality with $\alpha,\beta$ given by \eqref{eq:BK230-1}
implies that
\begin{align}
&
\left(
\frac{1}{m(J)}\sum_{k\in J}
\left(w_k^{\delta(\theta)/p(\theta)}\right)^{p(\theta)}
\right)^{1/p(\theta)}
=
\left(
\frac{1}{m(J)}\sum_{k\in J}
\left(w_k^{\delta_1}\right)^{1/\alpha}
\left(w_k^{\delta_2}\right)^{1/\beta}
\right)^{1/p(\theta)}
\nonumber\\
&\quad\le 
\left(
\left(\frac{1}{m(J)}\sum_{k\in J} w_k^{\delta_1}\right)^{1/\alpha}
\left(\frac{1}{m(J)}\sum_{k\in J} w_k^{\delta_2}\right)^{1/\beta}
\right)^{1/p(\theta)}
\nonumber\\
&\quad=
\left(
\frac{1}{m(J)}\sum_{k\in J} \left(w_k^{\delta_1/p_1}\right)^{p_1}
\right)^{\frac{1}{p_1}\frac{(1-\theta)p_1}{p(\theta)}}
\left(
\frac{1}{m(J)}\sum_{k\in J} \left(w_k^{\delta_2/p_2}\right)^{p_2}
\right)^{\frac{1}{p_2}\frac{\theta p_2}{p(\theta)}}.
\label{eq:BK230-2}
\end{align}
Define $q_j$ and $q(\theta)$ by
\[
\frac{1}{p_j}+\frac{1}{q_j}=1,
\quad
j=1,2,
\qquad
\frac{1}{p(\theta)}+\frac{1}{q(\theta)}=1,
\quad
\theta\in(0,1).
\]  
Let 
\begin{equation}\label{eq:BK230-3}
\gamma:=\frac{p(\theta)}{q(\theta)}\frac{q_1}{p_1}\frac{1}{1-\theta},
\quad
\mu:=\frac{p(\theta)}{q(\theta)}\frac{q_2}{p_2}\frac{1}{\theta}.
\end{equation}
Then
\begin{align*}
\frac{1}{\gamma}+\frac{1}{\mu}
&=
\frac{q(\theta)}{p(\theta)}
\left(\frac{p_1}{q_1}(1-\theta)+\frac{p_2}{q_2}\theta\right)
\\
&=
\frac{1}{p(\theta)-1}((p_1-1)(1-\theta)+(p_2-1)\theta))
=
\frac{p(\theta)-1}{p(\theta)-1}=1
\end{align*}
and
\begin{align*}
-\frac{\delta(\theta)q(\theta)}{p(\theta)}
&=
-\frac{q(\theta)}{p(\theta)}
((1-\theta)\delta_1+\theta\delta_2)
\\
&=
-\frac{q(\theta)}{p(\theta)}
\frac{p_1}{q_1}(1-\theta)
\delta_1 \frac{q_1}{p_1}
-\frac{q(\theta)}{p(\theta)}
\frac{p_2}{q_2}\theta
\delta_2 \frac{q_2}{p_2}
\\
&=
\left(-\frac{\delta_1q_1}{p_1}\right)\frac{1}{\gamma}
+
\left(-\frac{\delta_2 q_2}{p_2}\right)\frac{1}{\mu}.
\end{align*}
Hence, for all $k\in J$,
\[
\left(w_k^{-\delta(\theta)/p(\theta)}\right)^{q(\theta)}
=
\left(w_k^{-\delta_1q_1/p_1}\right)^{1/\gamma}
\left(w_k^{-\delta_2q_2/p_2}\right)^{1/\mu}.
\]
Applying H\"older's inequality with $\gamma$ and $\mu$ given by
\eqref{eq:BK230-3}, we get
\begin{align}
&
\left(
\frac{1}{m(J)}\sum_{k\in J}
\left(w_k^{-\delta(\theta)/p(\theta)}\right)^{q(\theta)}
\right)^{1/q(\theta)}
\nonumber\\
&\quad=
\left(
\frac{1}{m(J)}\sum_{k\in J}
\left(w_k^{-\delta_1q_1/p_1}\right)^{1/\gamma}
\left(w_k^{-\delta_2q_2/p_2}\right)^{1/\mu}
\right)^{1/q(\theta)}
\nonumber\\
&\quad\le 
\left(
\left(\frac{1}{m(J)}\sum_{k\in J}w_k^{-\delta_1q_1/p_1}\right)^{1/\gamma}
\left(\frac{1}{m(J)}\sum_{k\in J}w_k^{-\delta_2q_2/p_2}\right)^{1/\mu}
\right)^{1/q(\theta)}
\nonumber\\
&\quad=
\left(\frac{1}{m(J)}\sum_{k\in J}
\left(w_k^{-\delta_1/p_1}\right)^{q_1}\right)^{1/(\gamma q(\theta))}
\left(\frac{1}{m(J)}\sum_{k\in J}
\left(w_k^{-\delta_2/p_2}\right)^{q_2}\right)^{1/(\mu q(\theta))}
\nonumber\\
&\quad=
\left(\frac{1}{m(J)}\sum_{k\in J}
\left(w_k^{-\delta_1/p_1}\right)^{q_1}
\right)^{\frac{1}{q_1}\frac{(1-\theta)p_1}{p(\theta)}}
\left(\frac{1}{m(J)}\sum_{k\in J}
\left(w_k^{-\delta_2/p_2}\right)^{q_2}
\right)^{\frac{1}{q_2}\frac{\theta p_2}{p(\theta)}}.
\label{eq:BK230-4}
\end{align}
Multiplication of \eqref{eq:BK230-2} and \eqref{eq:BK230-4} shows that
\begin{align*}
&
\left(
\frac{1}{m(J)}\sum_{k\in J}
\left(w_k^{\delta(\theta)/p(\theta)}\right)^{p(\theta)}
\right)^{1/p(\theta)}
\left(
\frac{1}{m(J)}\sum_{k\in J}
\left(w_k^{-\delta(\theta)/p(\theta)}\right)^{q(\theta)}
\right)^{1/q(\theta)}
\\
&\quad\le 
\left(
\frac{1}{m(J)}\sum_{k\in J} \left(w_k^{\delta_1/p_1}\right)^{p_1}
\right)^{\frac{1}{p_1}\frac{(1-\theta)p_1}{p(\theta)}}
\left(\frac{1}{m(J)}\sum_{k\in J}
\left(w_k^{-\delta_1/p_1}\right)^{q_1}
\right)^{\frac{1}{q_1}\frac{(1-\theta)p_1}{p(\theta)}}
\\
&\qquad\times
\left(
\frac{1}{m(J)}\sum_{k\in J} \left(w_k^{\delta_2/p_2}\right)^{p_2}
\right)^{\frac{1}{p_2}\frac{\theta p_2}{p(\theta)}}
\left(\frac{1}{m(J)}\sum_{k\in J}
\left(w_k^{-\delta_2/p_2}\right)^{q_2}
\right)^{\frac{1}{q_2}\frac{\theta p_2}{p(\theta)}}
\\
&\quad\le 
\left(
\left[w^{\delta_1/p_1}\right]_{A_{p_1}(\mathbb{Z}_+)}
\right)^{\frac{(1-\theta)p_1}{p(\theta)}}
\left(
\left[w^{\delta_2/p_2}\right]_{A_{p_2}(\mathbb{Z}_+)}
\right)^{\frac{\theta p_2}{p(\theta)}}.
\end{align*}
Passing to the supremum over all intervals $J\subset\mathbb{Z}_+$,
we obtain
\[
\left[w^{\delta(\theta)/p(\theta)}\right]_{A_{p(\theta)}(\mathbb{Z}_+)}
\le 
\left(
\left[w^{\delta_1/p_1}\right]_{A_{p_1}(\mathbb{Z}_+)}
\right)^{\frac{(1-\theta)p_1}{p(\theta)}}
\left(
\left[w^{\delta_2/p_2}\right]_{A_{p_2}(\mathbb{Z}_+)}
\right)^{\frac{\theta p_2}{p(\theta)}}
\]
and $w^{\delta(\theta)/p(\theta)}\in A_{p(\theta)}(\mathbb{Z}_+)$.
\end{proof}
\subsection{Stability property of discrete Muckenhoupt weights}
The stability property of weights in $A_p(\mathbb{Z}_+)$ available in
\cite[Theorem~1.1]{BSey99} is not sufficient for our purposes. We will need
a more general version of it, which is analogous to \cite[Theorem~2.31]{BK97}.
\begin{theorem}\label{th:Muckenhoupt-stability}
Let $1<p_0<\infty$. If $w\in A_{p_0}(\mathbb{Z}_+)$, then there exists 
$\varepsilon_0>0$ such that $w^{1+\varepsilon}\in A_p(\mathbb{Z}_+)$ for all 
$\varepsilon\in(-\varepsilon_0,\varepsilon_0)$ and all 
$p\in(p_0-\varepsilon_0,p_0+\varepsilon_0)$.
\end{theorem}
\begin{proof}
The proof given below is an adaptation {of} the proof of 
\cite[Theorem~2.31]{BK97} to the discrete setting.
If $w\in A_{p_0}(\mathbb{Z}_+)$, then $w^{-1}\in A_{q_0}(\mathbb{Z}_+)$, where
$1/p_0+1/q_0=1$. By Theorem~\ref{th:RHI}, there are constants 
$\delta,\eta\in(0,\infty)$ and $C,D\in(0,\infty)$ such that
\begin{align}
\left(\frac{1}{m(R)}\sum_{k\in R}w_k^{p_0(1+\delta)}\right)^{1/(1+\delta)}
& \le 
\frac{C}{m(R)}\sum_{k\in R}w_k^{p_0},
\label{eq:Muckenhoupt-stability-1}
\\
\left(\frac{1}{m(R)}\sum_{k\in R}w_k^{-q_0(1+\eta)}\right)^{1/(1+\eta)}
& \le 
\frac{D}{m(R)}\sum_{k\in R}w_k^{-q_0}
\label{eq:Muckenhoupt-stability-2}
\end{align}
for all intervals $R\subset\mathbb{Z}_+$ 
{such that}
 $m(R)=2^r$ with
$r\in\mathbb{N}$.

If $-1<\lambda_1\le\lambda_2$, then $\Phi(x)=x^{(1+\lambda_2)/(1+\lambda_1)}$
is convex on $[0,\infty)$ and we may apply Jensen's inequality to deduce that
\[
\left(\frac{1}{m(R)}\sum_{k\in R}v_k^{1+\lambda_1}\right)^{1/(1+\lambda_1)}
\le 
\left(\frac{1}{m(R)}\sum_{k\in R}v_k^{1+\lambda_2}\right)^{1/(1+\lambda_2)}
\]
for any weight $v=\{v_k\}_{k\in\mathbb{Z}_+}$. Thus, letting
$\lambda:=\min\{\delta,\eta\}$, we obtain from 
\eqref{eq:Muckenhoupt-stability-1} and \eqref{eq:Muckenhoupt-stability-2}
that
\begin{align}
\left(\frac{1}{m(R)}\sum_{k\in R}w_k^{p_0(1+\lambda)}\right)^{1/(1+\lambda)}
&\le 
\left(\frac{1}{m(R)}\sum_{k\in R}w_k^{p_0(1+\delta)}\right)^{1/(1+\delta)}
\nonumber\\
&\le 
\frac{C}{m(R)}\sum_{k\in R}w_k^{p_0},
\label{eq:Muckenhoupt-stability-3}
\\
\left(\frac{1}{m(R)}\sum_{k\in R}w_k^{-q_0(1+\lambda)}\right)^{1/(1+\lambda)}
&\le
\left(\frac{1}{m(R)}\sum_{k\in R}w_k^{-q_0(1+\eta)}\right)^{1/(1+\eta)}
\nonumber\\
&\le 
\frac{D}{m(R)}\sum_{k\in R}w_k^{-q_0}
\label{eq:Muckenhoupt-stability-4}
\end{align}
for all intervals $R\subset\mathbb{Z}_+$ 
{such that}
 $m(R)=2^r$ with
$r\in\mathbb{N}$. It follows from \eqref{eq:Muckenhoupt-stability-3} and
\eqref{eq:Muckenhoupt-stability-4} that
\begin{align*}
&
\left(\frac{1}{m(R)}\sum_{k\in R}w_k^{p_0(1+\lambda)}\right)^{1/p_0}
\left(\frac{1}{m(R)}\sum_{k\in R}w_k^{-q_0(1+\lambda)}\right)^{1/q_0}
\\
&\quad\le
\left[
\left(\frac{C}{m(R)}\sum_{k\in R}w_k^{p_0}\right)^{1/p_0}
\left(\frac{D}{m(R)}\sum_{k\in R}w_k^{-q_0}\right)^{1/q_0}
\right]^{1+\lambda}
\\
&\quad \le
\left(C^{1/p_0}D^{1/q_0}\right)^{1+\lambda}[w]_{A_{p_0}(\mathbb{Z}_+)}^{1+\lambda}
\end{align*}
for all intervals $R\subset\mathbb{Z}_+$ 
{such that} $m(R)=2^r$ with $r\in\mathbb{N}$.

Now let $J\subset\mathbb{Z}_+$ be an arbitrary interval with $m(J)\ge 2$. 
Choose  any interval $R\subset\mathbb{Z}_+$ such that $m(R)=2^r$ for some 
$r\in\mathbb{N}$, $J\subset R$ and $m(R)\le 2m(J)$. Then
\begin{align*}
&
\frac{1}{m(J)}
\left(\sum_{k\in J}\left(w_k^{1+\lambda}\right)^{p_0}\right)^{1/p_0}
\left(\sum_{k\in J}\left(w_k^{1+\lambda}\right)^{-q_0}\right)^{1/q_0}
\\
&\quad\le
\frac{2}{m(R)}
\left(\sum_{k\in R}\left(w_k^{1+\lambda}\right)^{p_0}\right)^{1/p_0}
\left(\sum_{k\in R}\left(w_k^{1+\lambda}\right)^{-q_0}\right)^{1/q_0}
\\
&\quad\le
2\left(C^{1/p_0}D^{1/q_0}\right)^{1+\lambda}
[w]_{A_{p_0}(\mathbb{Z}_+)}^{1+\lambda}.
\end{align*}
Passing to the supermum over all intervals $J\subset\mathbb{Z}_+$, 
we obtain 
\[
\left[w^{1+\lambda}\right]_{A_{p_0}(\mathbb{Z}_+)}
\le 
2\left(C^{1/p_0}D^{1/q_0}\right)^{1+\lambda}
[w]_{A_{p_0}(\mathbb{Z}_+)}^{1+\lambda}.
\]
Consequently, $(p_0,(1+\lambda)p_0)\in\Gamma$. Because, 
by Lemma~\ref{le:BK230}, the set $\Gamma$ is convex and contains the 
half-line $(1,\infty)\times\{0\}$, it follows that $(p_0,p_0)$ is an 
inner point of $\Gamma$ (in the usual topology of $\mathbb{R}^2$), which 
is equivalent to the assertion of the theorem.
\end{proof}
Lemma~\ref{le:even-extension} and Theorem~\ref{th:Muckenhoupt-stability}
immediately imply the following.
\begin{corollary}\label{co:Muckenhoupt-stability}
Let $1<p_0<\infty$. If $w\in A_{p_0}^{\rm sym}(\mathbb{Z})$, then there exists 
$\varepsilon_0>0$ such that $w^{1+\varepsilon}\in A_p^{\rm sym}(\mathbb{Z})$ 
for all $\varepsilon\in(-\varepsilon_0,\varepsilon_0)$ and all 
$p\in(p_0-\varepsilon_0,p_0+\varepsilon_0)$.
\end{corollary}
\subsection{Stechkin's inequality}
In 1950, Stechkin \cite{S50} proved that if a $2\pi$-periodic function 
$a:\mathbb{R}\to\mathbb{Z}$ is  of bounded variation, then 
$a\in M_{\ell^p(\mathbb{Z})}$. B\"ottcher and  Seybold 
\cite[Theorem~1.2]{BSey99} proved the following weighted extension of 
Stechkin's inequality.
\begin{theorem}\label{th:Stechkin-weighted-ell-p}
Let $1<p<\infty$. If $w\in A_p(\mathbb{Z})$, then there is a constant
$c_{p,w}\in(0,\infty)$ depending only on $p$ and $w$ such that
\[
\|a\|_{M_{\ell^p(\mathbb{Z},w)}}
\le 
c_{p,w}\left(\|a\|_{L^\infty(-\pi,\pi)}+V(a)\right)
\] 
for all $a\in BV[-\pi,\pi]$. In particular, 
$BV[-\pi,\pi]\hookrightarrow M_{\ell^p(\mathbb{Z},w)}$.
\end{theorem}
\section{Banach algebra $M_{X(\mathbb{Z})}$}\label{sec:algebra-MX}
\subsection{Laurent operators on the space $\ell^2(\mathbb{Z})$}
We start this section with the following well-known result
(see, e.g., \cite[Theorem~1.1 and equalities (1.5)--(1.6)]{BS99} and also
\cite[Section~XXIII.2]{GGK93}, \cite[Section~3.1]{GGK03}).
\begin{theorem}\label{th:M2=L-infinity}
Let $\{a_n\}_{n\in\mathbb{Z}}$ be a sequence of complex numbers. The Laurent 
matrix $(a_{j-k})_{j,k=-\infty}^\infty$ generates a bounded operator $A$ on 
the space $\ell^2(\mathbb{Z})$ if and only if there exists a function 
$a\in L^\infty(-\pi,\pi)$ such that $a_n=\widehat{a}(n)$ for all 
$n\in\mathbb{Z}$. In this case $A=L(a)$ and
\[
\|a\|_{M_{\ell^2(\mathbb{Z})}}
=
\|L(a)\|_{\mathcal{B}(\ell^2(\mathbb{Z}))}
=
\|a\|_{L^\infty(-\pi,\pi)}.
\]
Moreover, if $a,b\in L^\infty(-\pi,\pi)$, then $L(a)L(b)=L(ab)$ on the space 
$\ell^2(\mathbb{Z})$. 
\end{theorem}
\subsection{Periodic distributions and their Fourier coefficients}
\label{sec:distributions}
Let $\mathcal{P}$ be the set of all infinitely differentiable $2\pi$-periodic
functions from $\mathbb{R}$ to $\mathbb{C}$. Elements of $\mathcal{P}$
are called periodic test functions. One can equip $\mathcal{P}$ with the 
countable family of seminorms
\[
\|u\|_{k,\mathcal{P}}:=\sup_{x\in[-\pi,\pi]}|D^{k-1}u(x)|,
\quad 
k\in\mathbb{N},
\]
where $D^ku$ denotes the $k$-th derivative of $u$ and $D^0u=u$, and the metric
\[
d(u,v):=\sum_{k=1}^\infty \frac{1}{2^k}
\frac{\|u-v\|_{k,\mathcal{P}}}{1+\|u-v\|_{k,\mathcal{P}}}.
\]
Then the set $\mathcal{P}$ endowed with the metric $d$ is a complete linear 
metric space (see \cite[Ch.~3, Theorems 2.1--2.2]{B73}). 

A periodic distribution is a continuous linear functional on the complete 
linear metric space $(\mathcal{P},d)$. The set of all periodic distributions 
is denoted by $\mathcal{P}'$. The functions 
\[
E_n(\theta):=e^{i n \theta},
\quad
\theta\in\mathbb{R},
\]
belong to $\mathcal{TP}\subset\mathcal{P}$ for all $n\in\mathbb{Z}$. 
The Fourier coefficients of a periodic distribution $a\in\mathcal{P}'$ 
are defined by
\[
\widehat{a}(n):=a(E_{-n}),
\quad
n\in\mathbb{Z}.
\]

The complex conjugate $\overline{a}$ of a periodic distribution 
$a\in\mathcal{P}'$ is defined by
\[
\overline{a}(u):=\overline{a(\overline{u})},
\quad
u\in\mathcal{P}.
\]
It is easy to see that 
\begin{equation}\label{eq:Fourier-coefficients-of-complex-conjugate}
\widehat{\overline{a}}(n)
=
\overline{a}(E_{-n})
=
\overline{a(\overline{E_{-n}})}
=
\overline{a(E_n)}
=
\overline{\widehat{a}(-n)},
\quad
n\in\mathbb{N}.
\end{equation}
\subsection{Complex conjugate of an element of $M_{X(\mathbb{Z})}$ belongs
to $M_{X'(\mathbb{Z})}$}
For $\varphi=\{\varphi_k\}_{k\in\mathbb{Z}}\in X(\mathbb{Z})$ and 
$\psi=\{\psi_k\}_{k\in\mathbb{Z}}\in X'(\mathbb{Z})$, put
\[
(\varphi,\psi):=\sum_{k\in\mathbb{Z}}\varphi_k \overline{\psi_k}.
\]
For $j\in\mathbb{Z}$, let $e_j:=\{\delta_{jk}\}_{k\in\mathbb{Z}}$, where
$\delta_{jk}$ is the Kronecker delta. It is clear that 
$e_k\in S_0(\mathbb{Z})$ for every $k\in\mathbb{Z}$.

The following lemma resembles \cite[Section~2.5(a)]{BS06}.
\begin{lemma}\label{le:multiplier-complex-conjugate}
Let $X(\mathbb{Z})$ be a Banach sequence space and $X'(\mathbb{Z})$ be its
associate space. If $a\in M_{X(\mathbb{Z})}$, then 
$\overline{a}\in M_{X'(\mathbb{Z})}$ and 
$\|a\|_{M_{X(\mathbb{Z})}}=\|\overline{a}\|_{M_{X'(\mathbb{Z})}}$.
\end{lemma}
\begin{proof}
If $i,j\in\mathbb{Z}$, then taking into account 
\eqref{eq:Fourier-coefficients-of-complex-conjugate}, we get
\begin{align*}
(a*e_i,e_j)
&=
\sum_{m\in\mathbb{Z}} (a*e_i)_m\overline{(e_j)_m}
=
(a*e_i)_j 
=
\sum_{k\in\mathbb{Z}}\widehat{a}(j-k)(e_i)_k
=\widehat{a}(j-i)
\\
&=
\overline{\widehat{\overline{a}}(i-j)}
=
\overline{\sum_{k\in\mathbb{Z}} \widehat{\overline{a}}(i-k)(e_j)_k}
=
\overline{(\overline{a}*e_j)_i}
=
\sum_{m\in\mathbb{Z}} (e_i)_m\overline{(\overline{a}*e_j)_m}
\\
&=
(e_i,\overline{a}*e_j).
\end{align*}
Hence, for all $\varphi,\psi\in S_0(\mathbb{Z})$,
\[
(a*\varphi,\psi)=(\varphi,\overline{a}*\psi).
\]
Fix $\psi\in S_0(\mathbb{Z})$ and recall that $X(\mathbb{Z})=X''(\mathbb{Z})$
with equal norms (see \cite[Ch.~1, Theorem~2.7]{BS88}). It follows 
from the above observations and Lemma~\ref{le:norm-on-nice-sequences} that
\begin{align}
\|\overline{a}*\psi\|_{X'(\mathbb{Z})}
&=
\sup\{|(\varphi,\overline{a}*\psi)|\ :\ 
\varphi\in S_0(\mathbb{Z}),\
\|\varphi\|_{X''(\mathbb{Z})}\le 1\}
\nonumber\\
&=
\sup\{|(a*\varphi,\psi)|\ :\ 
\varphi\in S_0(\mathbb{Z}),\
\|\varphi\|_{X(\mathbb{Z})}\le 1\}.
\label{eq:multiplier-complex-conjugate-1}
\end{align}
On the other hand, H\"older's inequality for Banach sequence spaces (see
\cite[Ch.~1, Theorem~2.4]{BS88}) and \eqref{eq:multiplier-norm} imply that
for all $\varphi\in S_0(\mathbb{Z})$,
\begin{equation}\label{eq:multiplier-complex-conjugate-2}
|(a*\varphi,\psi)|
\le
\|a*\varphi\|_{X(\mathbb{Z})}\|\psi\|_{X'(\mathbb{Z})}
\le 
\|a\|_{M_{X(\mathbb{Z})}}
\|\varphi\|_{X(\mathbb{Z})}
\|\psi\|_{X'(\mathbb{Z})}.
\end{equation}
Combining \eqref{eq:multiplier-complex-conjugate-1}
and \eqref{eq:multiplier-complex-conjugate-2}, we arrive at
\begin{equation}\label{eq:multiplier-complex-conjugate-3}
\|\overline{a}\|_{M_{X'(\mathbb{Z})}}
=
\sup\left\{
\frac{\|\overline{a}*\psi\|_{X'(\mathbb{Z})}}{\|\psi\|_{X'(\mathbb{Z})}}
\ :\
\psi\in S_0(\mathbb{Z}),\
\|\psi\|_{X'(\mathbb{Z})}\le 1
\right\}
\le 
\|a\|_{M_{X(\mathbb{Z})}}.
\end{equation}
Applying inequality \eqref{eq:multiplier-complex-conjugate-3}
to $\overline{a}$ and $X'(\mathbb{Z})$ in place of $a$ and 
$X(\mathbb{Z})$, respectively, we get
\begin{equation}\label{eq:multiplier-complex-conjugate-4}
\|a\|_{M_{X(\mathbb{Z})}}
=
\left\|\overline{\overline{a}}\right\|_{M_{X''(\mathbb{Z})}}
\le 
\|\overline{a}\|_{M_{X'(\mathbb{Z})}}.
\end{equation}
Inequalities \eqref{eq:multiplier-complex-conjugate-3} and
\eqref{eq:multiplier-complex-conjugate-4} imply that
$\|a\|_{M_{X(\mathbb{Z})}}=\|\overline{a}\|_{M_{X'(\mathbb{Z})}}$.
\end{proof}
\subsection{The sets $M_{X(\mathbb{Z})}$ and $M_{X'(\mathbb{Z})}$ coincide
if $X(\mathbb{Z})$ is reflection-invariant}
The next lemma extends \cite[Section~2.5(b)]{BS06}.
\begin{lemma}\label{le:multiplier-reflection-invariant-space}
If $X(\mathbb{Z})$ is a reflection-invariant Banach sequence space, then
$a\in M_{X(\mathbb{Z})}$ if and only if $\overline{a}\in M_{X(\mathbb{Z})}$
and $\|a\|_{M_{X(\mathbb{Z})}}=\|\overline{a}\|_{M_{X(\mathbb{Z})}}$.
\end{lemma}
\begin{proof}
If $X(\mathbb{Z})$ is reflection-invariant, then the anti-linear operator
$R$ defined by $R\varphi=\widetilde{\overline{\varphi}}$, that is,
\[
(R\varphi)_j=\overline{\varphi}_{-j},
\quad
j\in\mathbb{Z},
\]
is an isometry on $X(\mathbb{Z})$ and $R^2=I$.

If $\varphi\in S_0(\mathbb{Z})$, then taking into account 
\eqref{eq:Fourier-coefficients-of-complex-conjugate}, we get
for all $j\in\mathbb{Z}$,
\begin{align*}
(\overline{a}*\varphi)_j
&=
\sum_{k\in\mathbb{Z}}\widehat{\overline{a}}(j-k)\varphi_k
=
\sum_{k\in\mathbb{Z}}\overline{\widehat{a}(k-j)}\varphi_k
=
\sum_{k\in\mathbb{Z}}\overline{\widehat{a}(k-j)}\widetilde{\varphi}_{-k}
\\
&=
\overline{\sum_{k\in\mathbb{Z}}\widehat{a}(k-j)\widetilde{\overline{\varphi}}_{-k}}
=
\overline{\sum_{k\in\mathbb{Z}}\widehat{a}(-j-k)\widetilde{\overline{\varphi}}_k}
=
\overline{\sum_{k\in\mathbb{Z}}\widehat{a}(-j-k)(R\varphi)_k}
\\
&=
\overline{(a*R\varphi)}_{-j}
=
(R(a*R\varphi))_j.
\end{align*}
Hence, taking into account that $\varphi\in S_0(\mathbb{Z})$ if and only
if $\psi=R\varphi\in S_0(\mathbb{Z})$, we get
\begin{align*}
\|\overline{a}\|_{M_{X(\mathbb{Z})}}
&=
\sup\left\{
\frac{\|\overline{a}*\varphi\|_{X(\mathbb{Z})}}{\|\varphi\|_{X(\mathbb{Z})}}
\ :\ \varphi\in S_0(\mathbb{Z}),\ \varphi\ne 0
\right\}
\\
&=
\sup\left\{
\frac{\|R(a*R\varphi)\|_{X(\mathbb{Z})}}{\|\varphi\|_{X(\mathbb{Z})}}
\ :\ \varphi\in S_0(\mathbb{Z}),\ \varphi\ne 0
\right\}
\\
&=
\sup\left\{
\frac{\|a*R\varphi\|_{X(\mathbb{Z})}}{\|R(R\varphi)\|_{X(\mathbb{Z})}}
\ :\ \varphi\in S_0(\mathbb{Z}),\ \varphi\ne 0
\right\}
\\
&=
\sup\left\{
\frac{\|a*\psi\|_{X(\mathbb{Z})}}{\|R\psi\|_{X(\mathbb{Z})}}
\ :\ \psi\in S_0(\mathbb{Z}),\ \psi\ne 0
\right\}
\\
&=
\sup\left\{
\frac{\|a*\psi\|_{X(\mathbb{Z})}}{\|\psi\|_{X(\mathbb{Z})}}
\ :\ \psi\in S_0(\mathbb{Z}),\ \psi\ne 0
\right\}
=
\|a\|_{M_{X(\mathbb{Z})}},
\end{align*}
which completes the proof.
\end{proof}
The previous two lemmas lead to the following.
\begin{corollary}\label{co:multipliers-X-Xprime}
If $X(\mathbb{Z})$ is a reflection-invariant Banach sequence space and
$X'(\mathbb{Z})$ is its associate space, then
$M_{X(\mathbb{Z})}=M_{X'(\mathbb{Z})}$ and 
$\|a\|_{M_{X(\mathbb{Z})}}=\|a\|_{M_{X'(\mathbb{Z})}}$ for all
$a\in M_{X(\mathbb{Z})}$.
\end{corollary}
\begin{proof}
By Lemma~\ref{le:reflection-invariant-spaces}, $X'(\mathbb{Z})$
is a reflection-invariant Banach sequence space. If $a\in M_{X(\mathbb{Z})}$,
then $\overline{a}\in M_{X'(\mathbb{Z})}$ and 
$\|a\|_{M_{X(\mathbb{Z})}}=\|\overline{a}\|_{M_{X'(\mathbb{Z})}}$
in view of Lemma~\ref{le:multiplier-complex-conjugate}.
On the other hand, Lemma~\ref{le:multiplier-reflection-invariant-space} 
applied to $M_{X'(\mathbb{Z})}$ implies that $a\in M_{X'(\mathbb{Z})}$ and
$\|\overline{a}\|_{M_{X'(\mathbb{Z})}}
=
\|a\|_{M_{X'(\mathbb{Z})}}$. So,
\begin{equation}\label{eq:multipliers-X-Xprime-1}
M_{X(\mathbb{Z})}\hookrightarrow M_{X'(\mathbb{Z})}
\end{equation}
and
\begin{equation}\label{eq:multipliers-X-Xprime-2}
\|a\|_{M_{X(\mathbb{Z})}}
=
\|a\|_{M_{X'(\mathbb{Z})}}
\quad\mbox{for}\quad a\in M_{X(\mathbb{Z})}.
\end{equation}
It follows from what has already been proved with $X'(\mathbb{Z})$
in place of $X(\mathbb{Z})$ and the equality $X(\mathbb{Z})=X''(\mathbb{Z})$ 
with equal norms (see \cite[Ch.~1, Theorem~2.7]{BS88})
that
\begin{equation}\label{eq:multipliers-X-Xprime-3}
M_{X'(\mathbb{Z})}\hookrightarrow M_{X''(\mathbb{Z})}=M_{X(\mathbb{Z})}
\end{equation}
and
\begin{equation}\label{eq:multipliers-X-Xprime-4}
\|a\|_{M_{X'(\mathbb{Z})}}
=
\|a\|_{M_{X''(\mathbb{Z})}}
=
\|a\|_{M_{X(\mathbb{Z})}}
\quad\mbox{for}\quad a\in M_{X'(\mathbb{Z})}.
\end{equation}
Combining
\eqref{eq:multipliers-X-Xprime-1}--\eqref{eq:multipliers-X-Xprime-4},
we arrive at the desired conclusion.
\end{proof}
\subsection{Continuous embedding of $M_{X(\mathbb{Z})}$ into $L^\infty(-\pi,\pi)$}
The following result is the key ingredient in the proof of the main result
of this section. It extends one of the assertions in \cite[Section~2.5(d)]{BS06}.
\begin{theorem}\label{th:continuous-embedding}
Let $X(\mathbb{Z})$ be a reflexive reflection-invariant Banach sequence space.
If $a\in M_{X(\mathbb{Z})}$, then $a\in L^\infty(-\pi,\pi)$ and
\[
\|a\|_{L^\infty(-\pi,\pi)}\le \|a\|_{M_{X(\mathbb{Z})}}.
\]
\end{theorem}
\begin{proof}
If $a\in M_{X(\mathbb{Z})}$, then in view of 
Corollary~\ref{co:multipliers-X-Xprime}, $a\in M_{X'(\mathbb{Z})}$ and
$\|a\|_{M_{X(\mathbb{Z})}}=\|a\|_{M_{X'(\mathbb{Z})}}$. Since the space
$X(\mathbb{Z})$ is reflexive, Lemma~\ref{le:density}(b) implies that
$S_0(\mathbb{Z})$ is dense in $X(\mathbb{Z})$ and in $X'(\mathbb{Z})$.
Hence $L(a)\in\mathcal{B}(X(\mathbb{Z}))$ and $L(a)\in\mathcal{B}(X'(\mathbb{Z}))$
with
\[
\|L(a)\|_{\mathcal{B}(X(\mathbb{Z}))}
=
\|L(a)\|_{\mathcal{B}(X'(\mathbb{Z}))}
=
\|a\|_{M_{X(\mathbb{Z})}}
=
\|a\|_{M_{X'(\mathbb{Z})}}.
\]
It follows from Theorem~\ref{th:Calderon-Lozanovskii} that
\[
\|a\|_{M_{\ell^2(\mathbb{Z})}}
=
\|L(a)\|_{\mathcal{B}(\ell^2(\mathbb{Z}))}
\le 
\|L(a)\|_{\mathcal{B}(X(\mathbb{Z}))}^{1/2}
\|L(a)\|_{\mathcal{B}(X'(\mathbb{Z}))}^{1/2}
=
\|a\|_{M_{X(\mathbb{Z})}}.
\]
It remains to recall that $M_{\ell^2(\mathbb{Z})}=L^\infty(-\pi,\pi)$ and
$\|a\|_{M_{\ell^2(\mathbb{Z})}}=\|a\|_{L^\infty(-\pi,\pi)}$
(see Theorem~\ref{th:M2=L-infinity}).
\end{proof}
\subsection{$M_{X(\mathbb{Z})}$ is a Banach algebra if $X(\mathbb{Z})$
is reflexive and reflection-invariant}
Now we are in a position to prove the main result of this section.
It generalizes \cite[Section~2.5(g)]{BS06}.
\begin{theorem}\label{th:Banach-algebra}
Let $X(\mathbb{Z})$ be a reflexive reflection-invariant Banach sequence space.
Then $M_{X(\mathbb{Z})}$ is a Banach algebra under pointwise multiplication 
and the norm
\begin{equation}\label{eq:Banach-algebra-1}
\|a\|_{M_{X(\mathbb{Z})}}=\|L(a)\|_{\mathcal{B}(X(\mathbb{Z}))}.
\end{equation}
\end{theorem}
\begin{proof}
If $a,b\in M_{X(\mathbb{Z})}$, then by Theorem~\ref{th:continuous-embedding},
$a,b\in L^\infty(-\pi,\pi)$. Hence, Theorem~\ref{th:M2=L-infinity} implies
that
\[
L(ab)\varphi=L(a)L(b)\varphi
\]
for all $\varphi\in S_0(\mathbb{Z})\subset \ell^2(\mathbb{Z})$.
Since $S_0(\mathbb{Z})$ is dense in $X(\mathbb{Z})$ 
(see Lemma~\ref{le:density}(b)), we conclude that $L(ab)=L(a)L(b)$ on 
$X(\mathbb{Z})$. Hence $ab\in M_{X(\mathbb{Z})}$ and
\[
\|ab\|_{M_{X(\mathbb{Z})}}
=
\|L(ab)\|_{\mathcal{B}(X(\mathbb{Z}))}
\le
\|L(a)\|_{\mathcal{B}(X(\mathbb{Z}))}
\|L(b)\|_{\mathcal{B}(X(\mathbb{Z}))}
=
\|a\|_{M_{X(\mathbb{Z})}}
\|b\|_{M_{X(\mathbb{Z})}}.
\]
So, $M_{X(\mathbb{Z})}$ is a normed algebra under pointwise multiplication
and the norm \eqref{eq:Banach-algebra-1}.

It remains to show that $M_{X(\mathbb{Z})}$ is complete. Let 
$\{a^{(n)}\}_{n\in\mathbb{N}}$ be a Cauchy sequence in $M_{X(\mathbb{Z})}$.
Then Theorem~\ref{th:continuous-embedding} implies that 
$\{a^{(n)}\}_{n\in\mathbb{N}}$ is a Cauchy sequence in $L^\infty(-\pi,\pi)$.
Since the latter space is complete, there exists $a\in L^\infty(-\pi,\pi)$
such that $\|a^{(n)}-a\|_{L^\infty(-\pi,\pi)}\to 0$ as $n\to\infty$.
Theorem~\ref{th:M2=L-infinity} implies that
\begin{equation}\label{eq:Banach-algebra-2}
\|L(a^{(n)})-L(a)\|_{\mathcal{B}(\ell^2(\mathbb{Z}))}\to 0
\quad\mbox{as}\quad
n\to\infty.
\end{equation}
Since $\{L(a^{(n)})\}_{n\in\mathbb{N}}$ is a Cauchy sequence in 
$\mathcal{B}(X(\mathbb{Z}))$ and $\mathcal{B}(X(\mathbb{Z}))$ is complete, 
there exists an operator $A\in\mathcal{B}(X(\mathbb{Z}))$ such that
\begin{equation}\label{eq:Banach-algebra-3}
\|L(a^{(n)})-A\|_{\mathcal{B}(X(\mathbb{Z}))}\to 0
\quad\mbox{as}\quad
n\to\infty.
\end{equation}
Let $\varphi=\{\varphi_k\}_{k\in\mathbb{Z}}\in S_0(\mathbb{Z})
\subset\ell^2(\mathbb{Z})\cap X(\mathbb{Z})$. It follows from
\eqref{eq:Banach-algebra-2} and \eqref{eq:Banach-algebra-3} that 
\[
\|L(a^{(n)})\varphi-L(a)\varphi\|_{\ell^2(\mathbb{Z})}\to 0,
\quad
\|L(a^{(n)})\varphi-A\varphi\|_{X(\mathbb{Z})}\to 0
\quad\mbox{as}\quad
n\to\infty.
\]
In view of \cite[Ch.~1, Theorem~1.4]{BS88}, these relations imply that
for every finite sets $S\subset\mathbb{Z}$ and every $\varepsilon>0$,
\begin{align*}
&
m\big\{k\in S\ :\ 
\big|(L(a^{(n)})\varphi)_k-(L(a)\varphi)_k\big|\ge\varepsilon
\big\}
\to 0,
\\
&
m\big\{k\in S\ :\ 
\big|(L(a^{(n)})\varphi)_k-(A\varphi)_k\big|\ge\varepsilon
\big\}
\to 0
\end{align*}
as $n\to\infty$. Since $m(B)<1$ if and only if $B\subset\mathbb{Z}$ is empty,
we conclude that for every finite set $S\subset \mathbb{Z}$ and every
$\varepsilon>0$, there exists $n_0=n_0(S,\varepsilon)\in\mathbb{N}$
such that for all $k\in S$ and all $n>n_0$,
\[
\big|(L(a^{(n)})\varphi)_k-(L(a)\varphi)_k\big|<\varepsilon,
\quad
\big|(L(a^{(n)})\varphi)_k-(A\varphi)_k\big|<\varepsilon.
\]
Hence, for all $k\in\mathbb{Z}$, 
\[
(L(a)\varphi)_k
=
\lim_{n\to\infty} (L(a^{(n)}\varphi)_k
=
(A\varphi)_k.
\]
Thus $L(a)\varphi=A\varphi$ for all $\varphi\in S_0(\mathbb{Z})$. 
Lemma~\ref{le:density}(b) implies that $A=L(a)$ on $X(\mathbb{Z})$.
Since $L(a)\in\mathcal{B}(X(\mathbb{Z}))$, we conclude that 
$a\in M_{X(\mathbb{Z})}$. So, the Cauchy sequence 
$\{a^{(n)}\}_{n\in\mathbb{Z}}$ in $M_{X(\mathbb{Z})}$ converges to
$a\in M_{X(\mathbb{Z})}$. Thus $M_{X(\mathbb{Z})}$ is complete.
\end{proof}
\subsection{Uniform boundedness of the Fej\'er means in the norm of 
$M_{X(\mathbb{Z})}$}
Let $a\in\mathcal{P}'$ and let $\{\widehat{a}(n)\}_{n\in\mathbb{Z}}$ be the sequence
of Fourier coefficients of $a$. The partial sums of the Fourier series of $a$
are defined by
\[
[s_n(a)](\theta):=\sum_{k=-n}^n \widehat{a}(k)e^{ik\theta},
\quad
\theta\in\mathbb{R},
\quad
n\in\mathbb{Z}_+,
\]
and the Fej\'er (or Fej\'er-Ces\`aro) means of $a$ are defined by
\[
\sigma_n(a):=\frac{1}{n+1}\sum_{k=-n}^n s_k(a),
\quad n\in\mathbb{Z}_+.
\]

For $f,g\in L^1(-\pi,\pi)$, define the convolution of $f$ and $g$ by
\[
(f*g)(\theta):=\frac{1}{2\pi}\int_{-\pi}^\pi 
f(\theta-\phi)g(\phi)\,d\phi.
\]
For $n\in\mathbb{Z}_+$, let
\[
K_n(\theta)
:=
\sum_{|k|\le n}\left(1-\frac{|k|}{n+1}\right)e^{i\theta k}
=
\frac{1}{n+1}\left(
\frac{\sin\frac{n+1}{2}\theta}{\sin\frac{\theta}{2}}
\right)^2,
\quad
\theta\in\mathbb{R},
\]
be the $n$-th Fej\'er kernel. If $a\in L^1(-\pi,\pi)$, then
\begin{equation}\label{eq:Fejer-mean-as-convolution}
\sigma_n(a)=K_n* a
\end{equation}
(see, e.g., \cite[Ch.~I, Section 2]{K04}).

The following result for spaces $\ell^p(\mathbb{Z})$ is essentially due to 
Nikol'ski\u{\i} (see corollary to \cite[Theorem~5]{N66}). Its proof is given 
in \cite[Lemma~2.44]{BS06}. Our proof is analogous to the latter one.
\begin{lemma}\label{le:Fejer-uniform-boundedness}
Let $X(\mathbb{Z})$ be a reflexive reflection-invariant Banach sequence space. 
If $a\in M_{X(\mathbb{Z})}$, then
\[
\|\sigma_n(a)\|_{M_{X(\mathbb{Z})}}
\le 
\|a\|_{M_{X(\mathbb{Z})}}
\]
for all $n\in\mathbb{Z}_+$.
\end{lemma}
\begin{proof}
Fix $n\in\mathbb{Z}_+$. It follows from Theorem~\ref{th:continuous-embedding}
that $a\in L^\infty(-\pi,\pi)$. For $x\in\mathbb{R}$, define the 
function
\[
a_x(\theta):=a(\theta-x),
\quad
\theta\in\mathbb{R}.
\]
Then, for all $j\in\mathbb{Z}$,
\begin{align}
\widehat{a}_x(j)
&=
\frac{1}{2\pi}\int_{-\pi}^\pi a(\theta-x)e^{-ij\theta}\,d\theta
=
\frac{1}{2\pi}\int_{-\pi}^\pi a(\theta) e^{-ij(\theta+x)}\,d\theta
\nonumber\\
&=
e^{-ijx}\cdot\frac{1}{2\pi}\int_{-\pi}^\pi a(\theta)e^{-ij\theta}\,d\theta
=e^{-ijx}\widehat{a}(j).
\label{eq:Fejer-uniform-boundedness-1}
\end{align}
It follows from \eqref{eq:Fejer-mean-as-convolution} that
\begin{equation}\label{eq:Fejer-uniform-boundedness-2}
[\sigma_n(a)](\theta)
=
(K_n* a)(\theta)
=
\frac{1}{2\pi}\int_{-\pi}^\pi a_x(\theta)K_n(x) dx.
\end{equation}
For $x\in\mathbb{R}$, define the operator $D_x:X(\mathbb{Z})\to X(\mathbb{Z})$
by
\[
\{\varphi_j\}_{j\in\mathbb{Z}}\mapsto \{e^{ijx}\varphi_j\}_{j\in\mathbb{Z}}.
\]
It follows from the definition of the norm in $X(\mathbb{Z})$ that $D_x$
is an isometry on $X(\mathbb{Z})$.

If
$\varphi=\{\varphi_j\}_{j\in\mathbb{Z}}\in S_0(\mathbb{Z})$, then taking into 
account \eqref{eq:Fejer-uniform-boundedness-1}, we get
\begin{align}
(D_{-x} L(a)D_x)\varphi
&=
D_{-x}L(a)\left(\left\{e^{ijx}\varphi_j\right\}_{j\in\mathbb{Z}}\right)
=
D_{-x}\left(\left\{
\sum_{k\in\mathbb{Z}}\widehat{a}(j-k)e^{ikx}\varphi_k
\right\}_{j\in\mathbb{Z}}\right)
\nonumber\\
&=
\left\{
e^{-ijx}\sum_{k\in\mathbb{Z}}\widehat{a}(j-k)e^{ikx}\varphi_k
\right\}_{j\in\mathbb{Z}}
=
\left\{
\sum_{k\in\mathbb{Z}}\widehat{a}(j-k)e^{i(k-j)x}\varphi_k
\right\}_{j\in\mathbb{Z}}
\nonumber\\
&=
\left\{
\sum_{k\in\mathbb{Z}}\widehat{a}_x(j-k)\varphi_k
\right\}_{j\in\mathbb{Z}}
=
L(a_x)\varphi.
\label{eq:Fejer-uniform-boundedness-3}
\end{align}
It follows from the Lagrange mean value theorem that for all 
$\varphi=\{\varphi_j\}_{j\in\mathbb{Z}}\in S_0(\mathbb{Z})$, all 
$x,y\in\mathbb{R}$, and all $j\in\mathbb{Z}$,
\begin{align*}
|(D_x\varphi)_j-(D_y\varphi)_j|
&=
|(e^{ijx}-e^{ijy})\varphi_j|
\\
& 
{
=
\sqrt{[(\cos(jx)-\cos(jy)]^2+[\sin(jx)-\sin(jy)]^2}|\varphi_j|
}
\\
&\le 
{\sqrt{2}} |j|\,|x-y|\,|\varphi_j|
\le 
{\sqrt{2}}
\left(\max_{j\in\operatorname{supp}\varphi}|j|\right)|x-y|\,|\varphi_j|.
\end{align*}
The above inequality and Axioms (A1)--(A2) in the definition of a Banach
sequence space imply that
\begin{equation}\label{eq:Fejer-uniform-boundedness-4}
\|D_x\varphi-D_y\varphi\|_{X(\mathbb{Z})}
\le 
{\sqrt{2}}
\left(\max_{j\in\operatorname{supp}\varphi}|j|\right)
|x-y|\,\|\varphi\|_{X(\mathbb{Z})}.
\end{equation}
Equality \eqref{eq:Fejer-uniform-boundedness-3} and inequality
\eqref{eq:Fejer-uniform-boundedness-4} imply that for every
$\varphi\in S_0(\mathbb{Z})$, the function 
\begin{equation}\label{eq:Fejer-uniform-boundedness-5}
[-\pi,\pi]\to X(\mathbb{Z}),
\quad
x\mapsto L(a_x)\varphi
\end{equation}
is continuous. Then
\[
\frac{1}{2\pi}\int_{-\pi}^\pi \|L(a_x)\varphi\|_{X(\mathbb{Z})}K_n(x)\,dx
\le 
\max_{x\in[-\pi,\pi]}\|L(a_x)\varphi\|_{X(\mathbb{Z})}<\infty.
\]
So, by \cite[Ch.~2. Section~2, Theorem~2]{DU77}, the function given
in \eqref{eq:Fejer-uniform-boundedness-5} is Bochner integrable
on $[-\pi,\pi]$ with respect to the probability measure 
$d\mu_n(x)=K_n(x)\,dx/(2\pi)$. Similarly,
the constant function
\begin{equation}\label{eq:Fejer-uniform-boundedness-6}
[-\pi,\pi]\to X(\mathbb{Z}),
\quad
x\mapsto \varphi
\end{equation}
is continuous and, hence, is Bochner integrable. 

Equality \eqref{eq:Fejer-uniform-boundedness-2} and Fubini's theorem
imply that for all $j\in\mathbb{Z}$,
\begin{align*}
(\sigma_n(a))\widehat{\hspace{2mm}}(j)
&=
\frac{1}{2\pi}\int_{-\pi}^\pi [\sigma_n(a)](\theta)e^{-ij\theta}\,d\theta
\\
&=
\frac{1}{4\pi^2}\int_{-\pi}^\pi
\left(\int_{-\pi}^\pi a_x(\theta)K_n(x)\,dx\right)
e^{-ij\theta}\,d\theta
\\
&=
\frac{1}{4\pi^2}\int_{-\pi}^\pi
\left(\int_{-\pi}^\pi a_x(\theta)e^{-ij\theta}\,d\theta
\right)K_n(x)\,dx
\\
&=
\frac{1}{2\pi}\int_{-\pi}^\pi \widehat{a}_x(j)K_n(x)\,dx.
\end{align*}
Hence, for all $\varphi\in S_0(\mathbb{Z})$ and all $j\in\mathbb{Z}$,
\begin{align}
(L(\sigma_n(a))\varphi)_j
&=
\sum_{k\in\mathbb{Z}}(\sigma_n(a))\widehat{\hspace{2mm}}(j-k)\varphi_k
\nonumber\\
&=
\sum_{k\in\mathbb{Z}}
\left(\frac{1}{2\pi}\int_{-\pi}^\pi \widehat{a}_x(j-k)K_n(x)\,dx\right)
\varphi_k
\nonumber\\
&=
\frac{1}{2\pi}\int_{-\pi}^\pi \left(
\sum_{k\in\mathbb{Z}}\widehat{a}_x(j-k)\varphi_k
\right)K_n(x)\,dx
\nonumber\\
&=
\frac{1}{2\pi}\int_{-\pi}^\pi (L(a_x)\varphi)_j K_n(x)\,dx.
\label{eq:Fejer-uniform-boundedness-7}
\end{align}
Since the functions given by \eqref{eq:Fejer-uniform-boundedness-5}
and \eqref{eq:Fejer-uniform-boundedness-6} are Bochner integrable, it follows
from Hille's theorem (see 
\cite[Ch.~2, Section~2, {Theorem~6}]{DU77}) with
$T=T_j:X(\mathbb{Z})\to\mathbb{C}$, { sending a sequence
$\psi=\{\psi_i\}_{i\in\mathbb{Z}}\in X(\mathbb{Z})$ to its $j$-th coordinate}, that
\begin{align*}
\frac{1}{2\pi}\int_{-\pi}^\pi (L(a_x)\varphi)_j K_n(x)\,dx
&=
\int_{-\pi}^\pi T_j(L(a_x)\varphi)\,d\mu_n(x)
=
T_j\left(\int_{-\pi}^\pi L(a_x)\varphi\,d\mu_n(x)\right)
\\
&=
\left(\frac{1}{2\pi}\int_{-\pi}^\pi L(a_x)\varphi K_n(x)\,dx\right)_j.
\end{align*}
Combining this equality with \eqref{eq:Fejer-uniform-boundedness-7},
we get for all  $\varphi\in S_0(\mathbb{Z})$,
\begin{equation}\label{eq:Fejer-uniform-boundedness-8}
L(\sigma_n(a))\varphi
=
\frac{1}{2\pi}\int_{-\pi}^\pi L(a_x)\varphi K_n(x)\,dx.
\end{equation}

Let 
\[
f(j,x):=(L(a_x)\varphi)_j,
\quad
j\in\mathbb{Z},
\quad x\in[-\pi,\pi].
\]
It follows from \eqref{eq:Fejer-uniform-boundedness-8}, the Minkowski
integral inequality (see Lemma~\ref{le:Minkowski-inequality} with
the above function $f(j,x)$ and the probability measure 
$d\mu_n(x)=K_n(x)\,dx/(2\pi)$ on $[-\pi,\pi]$) that
\begin{align*}
\|L(\sigma_n(a))\varphi\|_{X(\mathbb{Z})}
&\le
\left\|\frac{1}{2\pi}\int_{-\pi}^\pi|L(a_x)\varphi| K_n(x)\,dx\right\|_{X(\mathbb{Z})}
\\
&\le 
\frac{1}{2\pi}\int_{-\pi}^\pi \|L(a_x)\varphi\|_{X(\mathbb{Z})}K_n(x)\,dx.
\end{align*}
Then applying \eqref{eq:Fejer-uniform-boundedness-3} and taking into account
that $D_{-x}$ and $D_x$ are isometries on $X(\mathbb{Z})$, we get
for all $\varphi\in S_0(\mathbb{Z})$,
\begin{align*}
\|L(\sigma_n(a))\varphi\|_{X(\mathbb{Z})}
&\le
\frac{1}{2\pi}\int_{-\pi}^\pi \|D_{-x}L(a)D_x\varphi\|_{X(\mathbb{Z})}K_n(x)\,dx
\\
&\le 
\|L(a)\|_{\mathcal{B}(X(\mathbb{Z}))}\|\varphi\|_{X(\mathbb{Z})}
\frac{1}{2\pi}\int_{-\pi}^\pi K_n(x)\,dx
\\
&=
\|a\|_{M_{X(\mathbb{Z})}}\|\varphi\|_{X(\mathbb{Z})}.
\end{align*}
Thus, 
\[
\|\sigma_n(a)\|_{M_{X(\mathbb{Z})}}
=
\sup\left\{
\frac{\|L(\sigma_n(a))\varphi\|_{X(\mathbb{Z})}}{\|\varphi\|_{X(\mathbb{Z})}}
:\varphi\in S_0(\mathbb{Z}),\varphi\ne 0\right\}
\le 
\|a\|_{M_{X(\mathbb{Z})}},
\]
which completes the proof.
\end{proof}
\section{Proof of the main results}\label{sec:proofs}
\subsection{Rearrangement-invariant Banach sequence spaces with symmetric weights}
Recall that a sequence $f=\{f_k\}_{k\in\mathbb{Z}}$ in a Banach sequence
space $X(\mathbb{Z})$ is said to have absolutely continuous norm in 
$X(\mathbb{Z})$ if $\|f\chi_{E_n}\|_{X(\mathbb{Z})}\to 0$ as $n\to\infty$
for every sequence $\{E_n\}_{n\in\mathbb{N}}$ of sets in $\mathbb{Z}$
satisfying $(\chi_{E_n})_k\to 0$ as $n\to\infty$ for all $k\in\mathbb{Z}$.
If all sequences $f=\{f_k\}_{k\in\mathbb{Z}}\in X(\mathbb{Z})$ have this 
property, then the space itself is said to have absolutely continuous norm
(see \cite[Ch.~1, Definition~3.1]{BS88}).
\begin{lemma}\label{le:ri-BSS-with-symmetric-weight}
If $X(\mathbb{Z})$ is a reflexive rearrangement-invariant Banach sequence
space and $w$ is a symmetric weight, then $X(\mathbb{Z},w)$ is a reflexive
reflection-invariant Banach sequence space.
\end{lemma}
\begin{proof}
Since $X(\mathbb{Z})$ is a reflexive Banach sequence, by 
\cite[Ch.~1, Corollary~4.4]{BS88}, $X(\mathbb{Z})$ and $X'(\mathbb{Z})$
have absolutely continuous norms. It is easy to see that then the
weighted Banach sequence spaces $X(\mathbb{Z},w)$ and $X'(\mathbb{Z},w^{-1})$
have absolutely continuous norms. Since $X'(\mathbb{Z},w^{-1})$ is the 
associate space of $X(\mathbb{Z},w)$, applying \cite[Ch.~1, Corollary~4.4]{BS88}
once again, we conclude that $X(\mathbb{Z},w)$ is reflexive. 

Let $f=\{f_k\}_{k\in\mathbb{Z}}\in X(\mathbb{Z})$. It is clear that
the counting measure $m$ on $\mathbb{Z}$ is reflection-invariant, that is,
$m(-S)=m(S)$ for every $S\subset\mathbb{Z}$. Since the weight $w$ is 
symmetric, we get for all $\lambda>0$,
\begin{align*}
d_{\widetilde{f}w}(\lambda)
&=
m\{k\in\mathbb{Z}:|\widetilde{f}_kw_k|>\lambda\}
=
m\{k\in\mathbb{Z}:|f_{-k}w_{-k}|>\lambda\}
\\
&=
m\{-k\in\mathbb{Z}:|f_{k}w_{k}|>\lambda\}
=
m\{k\in\mathbb{Z}:|f_{k}w_{k}|>\lambda\}
=d_{fw}(\lambda).
\end{align*}
Since $X(\mathbb{Z})$ is rearrangement-invariant, we obtain
\[
\|\widetilde{f}\|_{X(\mathbb{Z},w)}
=
\|\widetilde{f}w\|_{X(\mathbb{Z})}
=
\|fw\|_{X(\mathbb{Z})}
=
\|f\|_{X(\mathbb{Z},w)},
\]
that is, $X(\mathbb{Z},w)$ is reflection-invariant.
\end{proof}
\subsection{Proof of Theorem~\ref{th:Banach-algebra-of-multipliers}}
Part (a) follows from Lemma~\ref{le:ri-BSS-with-symmetric-weight} 
and Theorem~\ref{th:continuous-embedding}. Similarly, 
Lemma~\ref{le:ri-BSS-with-symmetric-weight} and Theorem~\ref{th:Banach-algebra}
yield part (b).

Let us prove part (c). Since $\alpha_X,\beta_X\in(0,1)$ and
$w\in A_{1/\alpha_X}^{\rm sym}(\mathbb{Z})\cap A_{1/\beta_X}^{\rm sym}(\mathbb{Z})$,
Corollary~\ref{co:Muckenhoupt-stability} yield that there exist $p$
and $q$ such that
\[
1<q<1/\beta_X\le 1/\alpha_X<p<\infty,
\quad
w\in A_p^{\rm sym}(\mathbb{Z})\cap A_q^{\rm sym}(\mathbb{Z}).
\]
By Theorem~\ref{th:Stechkin-weighted-ell-p}, there exist 
$c_{p,w},c_{q,w}\in(0,\infty)$ depending only on $p,q$ and $w$
such that for all $a\in BV[-\pi,\pi]$,
\begin{align}
\|L(a)\|_{\mathcal{B}(\ell^p(\mathbb{Z},w))}
&=
\|a\|_{M_{\ell^p(\mathbb{Z},w)}}
\le 
c_{p,w}\left(\|a\|_{L^\infty(-\pi,\pi)}+V(a)\right),
\label{eq:first-proof-1}
\\
\|L(a)\|_{\mathcal{B}(\ell^q(\mathbb{Z},w))}
&=
\|a\|_{M_{\ell^q(\mathbb{Z},w)}}
\le 
c_{q,w}\left(\|a\|_{L^\infty(-\pi,\pi)}+V(a)\right).
\label{eq:first-proof-2}
\end{align} 
On the other hand, Boyd's interpolation theorem (Theorem~\ref{th:Boyd})
implies that 
\begin{align}
\|a\|_{M_{X(\mathbb{Z},w)}}
&=
\|L(a)\|_{\mathcal{B}(X(\mathbb{Z},w))}
=
\|wL(a)w^{-1}\|_{\mathcal{B}(X(\mathbb{Z}))}
\nonumber\\
&\le 
C_{p,q}\max\left\{
\|wL(a)w^{-1}\|_{\mathcal{B}(\ell^p(\mathbb{Z}))},
\|wL(a)w^{-1}\|_{\mathcal{B}(\ell^q(\mathbb{Z}))}
\right\}
\nonumber\\
&=
C_{p,q}\max\left\{
\|L(a)\|_{\mathcal{B}(\ell^p(\mathbb{Z},w))},
\|L(a)\|_{\mathcal{B}(\ell^q(\mathbb{Z},w))}
\right\}.
\label{eq:first-proof-3}
\end{align}
Combining \eqref{eq:first-proof-1}--\eqref{eq:first-proof-3},
we arrive at
\[
\|a\|_{M_{X(\mathbb{Z},w)}}
\le 
C_{X(\mathbb{Z},w)}
\left(\|a\|_{L^\infty(-\pi,\pi)}+V(a)\right)
\]
with $C_{X(\mathbb{Z},w)}:=C_{p,q}c_{p,w}c_{q,w}$
independent of $a\in BV[-\pi,\pi]$.
\qed
\subsection{Convergence of Fej\'er means of continuous functions
of bounded variation in the multiplier norm}
Using interpolation techniques, we prove the following approximation
result, which plays a crucial role in the verification of property (iv)
in Theorem~\ref{th:Zalcman-Rudin} for $\mathcal{Y}=C_{X(\mathbb{Z},w)}$,
and, thus, in the proof of Theorem~\ref{th:algebra-C-plus-H-infinity}.

Let $C[-\pi,\pi]$ be the Banach space of all continuous $2\pi$-periodic
functions $f:\mathbb{R}\to\mathbb{C}$ with the norm
\[
\|f\|_{C[-\pi,\pi]}:=\sup_{\theta\in[-\pi,\pi]}|f(\theta)|.
\] 
\begin{lemma}\label{le:convergence-C-cap-BV-multiplier-norm}
Let $X(\mathbb{Z})$ be a reflexive rearrangement-invariant Banach sequence
space with nontrivial Boyd indices $\alpha_X,\beta_X$ and let
$w\in A_{1/\alpha_X}^{\rm sym}(\mathbb{Z})\cap 
A_{1/\beta_X}^{\rm sym}(\mathbb{Z})$. If $a\in C[-\pi,\pi]\cap BV[-\pi,\pi]$,
then
\[
\|\sigma_n(a)-a\|_{M_{X(\mathbb{Z},w)}}\to 0
\quad\mbox{as}\quad
n\to\infty.
\]
\end{lemma}
\begin{proof}
The idea of the proof is borrowed from \cite[Lemma~6.2]{BSey00}
and \cite[Theorems~3.2--3.3]{FKV21}. Fix $n\in\mathbb{Z}_+$.
Since $\alpha_X,\beta_X\in(0,1)$ and 
$w\in A_{1/\alpha_X}^{\rm sym}(\mathbb{Z})
\cap A_{1/\beta_X}^{\rm sym}(\mathbb{Z})$, it follows from 
Corollary~\ref{co:Muckenhoupt-stability} that there exist $p$ and $q$ such 
that
\[
1<q<1/\beta_X \le 1/\alpha_X<p<\infty,
\quad
w\in A_p^{\rm sym}(\mathbb{Z})\cap A_q^{\rm sym}(\mathbb{Z}).
\]
It follows from Boyd's interpolation theorem (see Theorem~\ref{th:Boyd})
that 
\begin{align}
&
\|\sigma_n(a)-a\|_{M_{X(\mathbb{Z},w)}}
=
\|L(\sigma_n(a)-a)\|_{\mathcal{B}(X(\mathbb{Z},w))}
=
\|wL(\sigma_n(a)-a)w^{-1}\|_{\mathcal{B}(X(\mathbb{Z}))}
\nonumber\\
&\quad\le 
C_{p,q}\max\left\{
\|wL(\sigma_n(a)-a)w^{-1}\|_{\mathcal{B}(\ell^p(\mathbb{Z}))},
\|wL(\sigma_n(a)-a)w^{-1}\|_{\mathcal{B}(\ell^q(\mathbb{Z}))}
\right\}
\nonumber\\
&\quad=
C_{p,q}\max\left\{
\|L(\sigma_n(a)-a)\|_{\mathcal{B}(\ell^p(\mathbb{Z},w))},
\|L(\sigma_n(a)-a)\|_{\mathcal{B}(\ell^q(\mathbb{Z},w))}
\right\}.
\label{eq:convergence-C-cap-BV-multiplier-norm-1}
\end{align}
Since $w\in A_p^{\rm sym}(\mathbb{Z})$, applying 
Corollary~\ref{co:Muckenhoupt-stability} once again, we conclude that 
$w^{1+\delta_2}\in A_{p(1+\delta_1)}^{\rm sym}(\mathbb{Z})$ whenever 
$|\delta_1|$ and $|\delta_2|$ are sufficiently small. If $p\ge 2$,
then one can find sufficiently small $\delta_1,\delta_2$ and 
{a number} $\theta_p\in(0,1)$ such that
\begin{equation}\label{eq:convergence-C-cap-BV-multiplier-norm-2}
\frac{1}{p}=\frac{1-\theta_p}{2}+\frac{\theta_p}{p(1+\delta_1)},
\quad
w=1^{1-\theta_p}w^{(1+\delta_2)\theta_p},
\quad
w^{1+\delta_2}\in A_{p(1+\delta_1)}^{\rm sym}(\mathbb{Z}).
\end{equation}
If $1<p<2$, then one can find a sufficiently small number $\delta_2>0$, 
a number $\delta_1<0$ with a sufficiently small $|\delta_1|$, and a 
number $\theta_p\in(0,1)$ such that all conditions in 
\eqref{eq:convergence-C-cap-BV-multiplier-norm-2} are fulfilled.
Taking into account \eqref{eq:convergence-C-cap-BV-multiplier-norm-2},
we obtain from the Stein-Weiss interpolation theorem (see, e.g.,
\cite[Ch.~4, Theorem~3.6]{BS88}) that
\begin{align}
&
\|L(\sigma_n(a)-a)\|_{\mathcal{B}(\ell^p(\mathbb{Z},w))}
\nonumber\\
&\quad\le 
\|
L(\sigma_n(a)-a)
\|_{\mathcal{B}(\ell^2(\mathbb{Z}))}^{1-\theta_p}
\|
L(\sigma_n(a)-a)
\|_{\mathcal{B}(\ell^{p(1+\delta_1)}(\mathbb{Z},w^{1+\delta_2}))}^{\theta_p}.
\label{eq:convergence-C-cap-BV-multiplier-norm-3}
\end{align}
It follows from Theorem~\ref{th:M2=L-infinity} that
\begin{equation}\label{eq:convergence-C-cap-BV-multiplier-norm-4}
\|L(\sigma_n(a)-a)\|_{\mathcal{B}(\ell^2(\mathbb{Z}))}
=
\|\sigma_n(a)-a\|_{L^\infty(-\pi,\pi)}
\le 
\|\sigma_n(a)-a\|_{C[-\pi,\pi]}.
\end{equation}
On the other hand, Lemmas~\ref{le:ri-BSS-with-symmetric-weight}
and~\ref{le:Fejer-uniform-boundedness} imply that
\begin{align}
\|
L(\sigma_n(a)-a)
\|_{\mathcal{B}(\ell^{p(1+\delta_1)}(\mathbb{Z},w^{1+\delta_2}))}
&=
\|\sigma_n(a)-a\|_{M_{\ell^{p(1+\delta_1)}(\mathbb{Z},w^{1+\delta_2})}}
\nonumber\\
&\le 
2\|a\|_{M_{\ell^{p(1+\delta_1)}(\mathbb{Z},w^{1+\delta_2})}}.
\label{eq:convergence-C-cap-BV-multiplier-norm-5}
\end{align}
Theorem~\ref{th:Stechkin-weighted-ell-p} yields that there exists
a constant $C(p,w,\delta_1,\delta_2)>0$ such that
\begin{equation}\label{eq:convergence-C-cap-BV-multiplier-norm-6}
\|a\|_{M_{\ell^{p(1+\delta_1)}(\mathbb{Z},w^{1+\delta_2})}}
\le 
C(p,w,\delta_1,\delta_2)
\left(\|a\|_{L^\infty(-\pi,\pi)}+V(a)\right).
\end{equation}
Combining \eqref{eq:convergence-C-cap-BV-multiplier-norm-3}%
--\eqref{eq:convergence-C-cap-BV-multiplier-norm-6}, we get 
\begin{align}
&
\|L(\sigma_n(a)-a))\|_{\mathcal{B}(\ell^p(\mathbb{Z},w))}
\nonumber\\
&\quad\le 
(2C(p,w,\delta_1,\delta_2))^{\theta_p} 
\|\sigma_n(a)-a\|_{C[-\pi,\pi]}^{1-\theta_p}
\left(
\|a\|_{L^\infty(-\pi,\pi)}+V(a)
\right)^{\theta_p}.
\label{eq:convergence-C-cap-BV-multiplier-norm-7}
\end{align}
Analogously, one can show that  there exist $\delta_3,\delta_4\in\mathbb{R}$
and $\theta_q\in(0,1)$ such that $|\delta_3|, |\delta_4|$ are sufficiently 
small and
\[
\frac{1}{q}=\frac{1-\theta_q}{2}+\frac{\theta_q}{q(1+\delta_3)},
\quad
w=1^{1-\theta_q}w^{(1+\delta_4)\theta_q},
\quad
w^{1+\delta_4}\in A_{q(1+\delta_3)}^{\rm sym}(\mathbb{Z}).
\]
Hence
\begin{align}
&
\|L(\sigma_n(a)-a))\|_{\mathcal{B}(\ell^q(\mathbb{Z},w))}
\nonumber\\
&\quad\le 
(2C(q,w,\delta_3,\delta_4))^{\theta_q} 
\|\sigma_n(a)-a\|_{C[-\pi,\pi]}^{1-\theta_q}
\left(
\|a\|_{L^\infty(-\pi,\pi)}+V(a)
\right)^{\theta_q},
\label{eq:convergence-C-cap-BV-multiplier-norm-8}
\end{align}
where $C(q,w,\delta_3,\delta_4)>0$ is some constant depending only on
$q,w$ and $\delta_3,\delta_4$.
It follows from \eqref{eq:convergence-C-cap-BV-multiplier-norm-1} and
\eqref{eq:convergence-C-cap-BV-multiplier-norm-7}--%
\eqref{eq:convergence-C-cap-BV-multiplier-norm-8} that for all $n\in\mathbb{Z}_+$,
\begin{align}
&
\|\sigma_n(a)-a\|_{M_{X(\mathbb{Z},w)}}
\nonumber\\
&\quad\le 
C\max\left\{
\left(
\|a\|_{L^\infty(-\pi,\pi)}+V(a)
\right)^{\theta_p},
\left(
\|a\|_{L^\infty(-\pi,\pi)}+V(a)
\right)^{\theta_q}
\right\}
\nonumber\\
&\qquad\times
\max\left\{
\|\sigma_n(a)-a\|_{C[-\pi,\pi]}^{1-\theta_p},
\|\sigma_n(a)-a\|_{C[-\pi,\pi]}^{1-\theta_q}
\right\},
\label{eq:convergence-C-cap-BV-multiplier-norm-9}
\end{align}
where
\[
C:=C_{p,q}\max\left\{
(2C(p,w,\delta_1,\delta_2))^{\theta_p},
(2C(q,w,\delta_3,\delta_4))^{\theta_q}
\right\}.
\]
It is well known that for $a\in C[-\pi,\pi]$,
\begin{equation}\label{eq:convergence-C-cap-BV-multiplier-norm-10}
\|\sigma_n(a)-a\|_{C[-\pi,\pi]}\to 0
\quad\mbox{as}\quad
n\to\infty
\end{equation}
(see, e.g., \cite[Ch.~I, Subsections~2.5, 2.10, 2.11]{K04}).
Combining \eqref{eq:convergence-C-cap-BV-multiplier-norm-9}
and
\eqref{eq:convergence-C-cap-BV-multiplier-norm-10}, we arrive at the desired 
conclusion.
\end{proof}
\subsection{Proof of Theorem~\ref{th:algebra-C-plus-H-infinity}}
The idea of the proof is borrowed from \cite[Theorem~2.53]{BS06}.
It follows from Theorem~\ref{th:Banach-algebra-of-multipliers}
that $M_{X(\mathbb{Z},w)}$ is a Banach algebra. We apply the Zalcman-Rudin 
theorem (see Theorem~\ref{th:Zalcman-Rudin}) to 
\[
\mathcal{X}=M_{X(\mathbb{Z},w)},
\quad
\mathcal{Y}=C_{X(\mathbb{Z},w)}, 
\quad
\mathcal{Z}=H_{X(\mathbb{Z},w)}^{\infty,\pm}, 
\]
and
\[
\Phi=\left\{S_n\in\mathcal{B}(M_{X(\mathbb{Z},w)})\ :\
S_na:=\sigma_n(a),\ n\in\mathbb{Z}_+\right\},
\]
where $\sigma_n(a)$ denotes the $n$-th Fej\'er mean of $a$. It is clear that
properties (i) and (ii) in Theorem~\ref{th:Zalcman-Rudin} are satisfied.
Property (iii) follows from Lemmas~\ref{le:ri-BSS-with-symmetric-weight} 
and~\ref{le:Fejer-uniform-boundedness}. It remains to verify property (iv).
It follows from the definition of $C_{X(\mathbb{Z},w)}$ that given 
$a\in C_{X(\mathbb{Z},w)}$ and $\varepsilon>0$, there exists $f\in\mathcal{TP}$
such that $\|a-f\|_{M_{X(\mathbb{Z},w)}}<\varepsilon/3$. Then
Lemmas~\ref{le:ri-BSS-with-symmetric-weight} 
and~\ref{le:Fejer-uniform-boundedness} imply that
\begin{align}
\|\sigma_n(a)-a\|_{M_{X(\mathbb{Z},w)}}
&\le 
\|\sigma_n(a-f)\|_{M_{X(\mathbb{Z},w)}}
+
\|\sigma_n(f)-f\|_{M_{X(\mathbb{Z},w)}}
+
\|f-a\|_{M_{X(\mathbb{Z},w)}}
\nonumber\\
& \le
2\|a-f\|_{M_{X(\mathbb{Z},w)}}
+
\|\sigma_n(f)-f\|_{M_{X(\mathbb{Z},w)}}
\nonumber\\
&\le 
\frac{2\varepsilon}{3}+\|\sigma_n(f)-f\|_{M_{X(\mathbb{Z},w)}}.
\label{eq:algebra-C-plus-H-infinity-1}
\end{align}
On the other hand, since $f\in\mathcal{TP}\subset C[-\pi,\pi]\cap BV[-\pi,\pi]$,
it follows from Lemma~\ref{le:convergence-C-cap-BV-multiplier-norm} that
\begin{equation}\label{eq:algebra-C-plus-H-infinity-2}
\|\sigma_n(f)-f\|_{M_{X(\mathbb{Z},w)}}<\frac{\varepsilon}{3}
\end{equation}
whenever $n$ is large enough. Combining \eqref{eq:algebra-C-plus-H-infinity-1}
and \eqref{eq:algebra-C-plus-H-infinity-2}, we see that
\[
\|S_na-a\|_{M_{X(\mathbb{Z},w)}}
=
\|\sigma_n(a)-a\|_{M_{X(\mathbb{Z},w)}}
<
\varepsilon
\]
whenever $n$ is large enough. So, property (iv)
in Theorem~\ref{th:Zalcman-Rudin} is satisfied. Thus, that theorem implies
that $C_{X(\mathbb{Z},w)}+H_{X(\mathbb{Z},w)}^{\infty,\pm}$ 
are closed subspaces of the Banach algebra $M_{X(\mathbb{Z},w)}$.

Now let $a,b\in C_{X(\mathbb{Z},w)}+H_{X(\mathbb{Z},w)}^{\infty,\pm}$. Then
there are $a_k,b_k\in\mathcal{TP}+H_{X(\mathbb{Z},\pm)}^{\infty,\pm}$ for 
$k\in\mathbb{N}$ such that
\begin{equation}\label{eq:algebra-C-plus-H-infinity-3}
\|a-a_k\|_{M_{X(\mathbb{Z},w)}}\to 0,
\quad
\|b-b_k\|_{M_{X(\mathbb{Z},w)}}\to 0
\quad\mbox{as}\quad
k\to\infty.
\end{equation}
It follows from Lemma~\ref{le:ri-BSS-with-symmetric-weight} and
Theorem~\ref{th:continuous-embedding} that 
$a_k,b_k\in L^\infty(-\pi,\pi)\subset L^2(-\pi,\pi)$. So $a_k$ and $b_k$
can be represented by their Fourier series convergent in $L^2(-\pi,\pi)$.
Looking at the Fourier series representation of $a_kb_k\in L^\infty(-\pi,\pi)$,
we conclude that $a_kb_k\in\mathcal{TP}+H_{X(\mathbb{Z},w)}^{\infty,\pm}$ for all
$k\in\mathbb{N}$. It follows from \eqref{eq:algebra-C-plus-H-infinity-3}
that 
\[
\|a_kb_k-ab\|_{M_{X(\mathbb{Z},w)}} \to 0
\quad\mbox{as}\quad
k\to\infty.
\]
Since $C_{X(\mathbb{Z},w)}+H_{X(\mathbb{Z},w)}^{\infty,\pm}$ is a closed
subspace of the Banach algebra $M_{X(\mathbb{Z},w)}$, we conclude that 
$ab\in C_{X(\mathbb{Z},w)}+H_{X(\mathbb{Z},w)}^{\infty,\pm}$.
Thus, $C_{X(\mathbb{Z},w)}+H_{X(\mathbb{Z},w)}^{\infty,\pm}$ is an algebra.
\qed
\subsection*{Acknowledgments}
We would like to thank the anonymous referees for useful remarks.
\subsection*{Funding}
This work is funded by national funds through the FCT - Funda\c{c}\~ao para a 
Ci\^encia e a Tecnologia, I.P., under the scope of the projects 
UIDB/00297/2020 (\url{https://doi.org/10.54499/UIDB/00297/2020})
and 
UIDP/ 00297/2020
(\url{https://doi.org/10.54499/UIDP/00297/2020})
(Center for Mathematics and Applications).
The second author is funded by national funds through the FCT – 
Funda\c{c}\~ao para a Ci\^encia e a Tecnologia, I.P., 
under the scope of the PhD scholarship UI/BD/ 152570/2022.
\bibliographystyle{abbrv}
\bibliography{OKST}
\end{document}